\documentstyle[12pt]{article} 
\begin{document} 
\title{Generalized Loop Groups of Complex
Manifolds, Gaussian Quasi-Invariant Measures on them
and their Representations.}  
\author{S.V. Ludkovsky.}  
\date{25 September 1999}
\maketitle
\par Istituto di Matematica dell'Universit\`a di Trieste,
\par Piazzale Europa, 1
\par I-34100, Trieste, Italia
\par Permament addresses: 
\par Theoretical Department, Institute of General Physics,
\par Str. Vavilov 38, Moscow, 117942, Russia. \\
\par 1980 Mathematics Subject Classification (1991 Revision): 
22A10, 43A05 and 43A65
\begin{abstract} 
\par Loop groups $G$ as families of
mappings of the complex manifold $M$
into another complex manifold $N$
preserving marked points $s_0\in M$ and $y_0\in N$
are investigated.  Quasi-invariant measures $\mu $ on $G$
relative to dense subgroups $G'$ are
constructed.  These measures are used for the
studying of irreducible representations
of such groups.
\end{abstract} 
\section{Introduction.}  
\par Loop groups are
very important in differential geometry, algebraic topology and theoretical
physics \cite{bryl,chen,mensk,swit}, but about Gaussian
quasi-invariant differentiable measures on them
nothing was known.  Only the 
simplest possible representations associated with
path's integrals were constructed for loop groups
of the circle, that is, for the manifold 
$M=S^1$ and (real) Riemann manifolds $N$ \cite{mensk}.  
On the other hand, the
quasi-invariant measures may be used for a construction of regular unitary
representations \cite{luum96,luum98,lutm99,lurim2}.  
Moreover, the quasi-invariant measures are helpful for an
investigation of the group itself.
In the previous papers of the author \cite{lulgihp,ludan}
loop groups of Riemann manifolds $M$ and $N$
were investigated, where either
$M=S^n$ were $n$-dimensional real spheres, $n=1,2,...$,
or $M=S^{\infty }$ was a unit sphere in a real separable
Hilbert space $l_2({\bf R})$. 
It was a progress in comparison
with previous works of others authors, which considered only loop groups
for the simplest case $M=S^1$.
\par  This article treats 
arbitrary complex separable connected metrizable manifolds $M$ and $N$.
For example, products of odd-dimensional
real spheres $S^{2n-1}\times S^{2m-1}$ may be supplied in different ways
by  structures of complex manifolds \cite{lnma}.
Another numerous examples of complex manifolds may be found in \cite{kob}
and references therein such as domains in $\bf C^n$, a complex torus
${\bf C^n}/D$, where $D$ is a discrete additive subgroup
of $\bf C^n$ generated by a basis $\tau _1,...,\tau _{2n}$ of 
$\bf C^n$ over $\bf R$; a quotient space $G/D$ of a complex Lie group
$G$ by a discrete subgroup $D$, submanifolds of the complex
Grassman manifold $G_{p,q}({\bf C}),$ also their different products
and their submanifolds.
In general there are complex compact manifolds, which are not 
K\"ahler manifolds \cite{loni,meve}.
For the construction of loop groups here are used manifolds $M$
with some mild additional conditions.
When $M$ is finite-dimensional over $\bf C$ we suppose that it is compact.
This condition is not very restrictive, since each locally compact
space has Alexandroff (one-point) compactification (see Theorem 3.5.11
in \cite{eng}).
When $M$ is infinite-dimensional over $\bf C$ it is assumed,
that $M$ is embedded as a closed bounded subset into the 
corresponding Banach space $X_M$ over $\bf C$.
This is necessary that to define a group structure on a 
quotient space of a free loop space. 
\par The free loop space
is considered as consisting of continuous functions $f: M\to N$
which are holomorphic on $M\setminus M'$ and preserving marked points
$f(s_0)=y_0$, where $M'$
is a closed real submanifold depending on $f$ with a codimension
$codim_{\bf R}M'=1$, $s_0\in M$ and $y_0\in N$ are marked points. 
There are two reasons to consider
such class of mappings. The first is the need to define correctly
compositions of elements in the loop group (see beneath).
The second is the isoperimetric inequality for holomorphic
loops, which can cause the condition of a loop
to be constant on a sufficiently small neighbourhood of $s_0$
in $M$, if this loop is in some small neighbourhood of $w_0$,
where $w_0(M):=\{ y_0 \} $ is a constant loop 
(see Remark 3.2 in \cite{hum}).
\par In this article loop groups of different classes are
considered. Classes analogous to Gevrey classes
of $f: M\setminus M' \to N$  
are considered for the construction of dense loop 
subgroups and quasi-invariant measures. Henceforth, we consider
only orientable manifolds $M$ and $N$, since 
for a non-orientable manifold there always exists
its orientable double covering manifold (see \S 6.5 in \cite{abma}).
Loop commutative monoids with the cancellation property
are quotients of families of mappings $f$ from $M$ into a
manifold $N$ with $f(s_0)=y_0$ by the corresponding equivalence
relation. 
For the definition of the equivalence relation here are 
not used groups of holomorphic diffeomorphisms because of strong
restrictions on their structure caused by 
holomorphicity (see Theorems 1 and 2 in \cite{bomon}).
Groups are constructed from monoids with the help of A. Grothendieck
procedure.
These groups are commutative and non-locally compact.
They does not have non-trivial local one-parameter subgroups
$\{ g^b:$ $b\in (-a,a) \} $ with $a>0$ for an element $g$ corresponding
to a class of a mapping $f: M\to N$, $f(s_0)=y_0$, when 
$f$ is such that $\sup_{y\in N} [card(f^{-1}(y))]=k<\aleph _0$, since 
$g^{1/p}$ does not exist in the loop group for each prime integer
$p$ such that $p>k$ (see \S 2).
Therefore, they are not Banach-Lie groups,
since in each neighbourhood
$W$ of the unit element $e$ there are elements which does not belong
to any local one-parameter subgroup.
\par Irreducible components of strongly continuous unitary representations 
of Abelian locally compact groups are one-dimensional by Theorem 22.17
\cite{hew}. In general commutative non-locally compact groups may 
have infinite-dimensional irreducible strongly continuous unitary 
representations, for example, infinite-dimensional Banach spaces 
over $\bf R$ considered as additive groups
(see \S 2.4 in \cite{bana} and \S 4.5 \cite{gelv} ).
\par Quasi-invariant measures are constructed on these loop groups. 
Then such measures
are used for the investigation of irreducible unitary representations.  
Loop groups, their structure and quasi-invariant measures on them
are investigated in \S 2.  
Unitary representations of loop groups are given in \S 3.
Irreducibility of unitary representations of a dense subgroup 
$G'$ associated with purely Gaussian quasi-invariant 
infinitely differentiable measures on the entire group
$G$ is investigated using specific properties of a quasi-invariance factor
relative to shifts from the dense subgroup $G'$.
Characters and also infinite-dimensional irreducible 
unitary representations are investigated below.
The relation between an equivalence of regular representations 
and equivalence of the corresponding measures is studied.
\section{Loop groups.}  
\par To avoid misunderstanding we
first give our Definitions and notations.  
\par {\bf 2.1. Definitions and Notes.1.}  Let $G$ be a 
Hausdorff topological group, we denote by 
$\mu : Af(G,\mu )\to [0,\infty )\subset \bf R$ 
a $\sigma $-additive measure. Its left shifts 
$\mu _{\phi }(E):=\mu (\phi ^{-1}\circ E)$ are considered for each 
$E \in Af(G,\mu )$,
where $Af(G,\mu )$ is
the completion of $Bf(G)$ by $\mu $-null sets, 
$Bf(G)$ is the Borel $\sigma $-field on $G$,
$\phi \circ E:=\{ \phi \circ h: h\in E \} $.
Then $\mu $ is called quasi-invariant if there exists a dense subgroup
$G'$ such that $\mu _{\phi }$ is equivalent to $\mu $ for each $\phi \in G'$.
Henceforth, we assume that a
quasi-invariance factor $\rho _{\mu }(\phi ,g)=\mu _{\phi }(dg)/\mu (dg)$
is continuous by $(\phi ,g) \in G' \times G$,
$\mu (V)>0$ for some (open) neighbourhood $V\subset
G$ of the unit element $e \in G$ and
$\mu (G)<\infty $. 
\par Let $({\sf M,F})$ be a space $\sf M$ of measures on $(G,Bf(G))$
with values
in $\bf R$ and $G"$ be a dense subgroup in $G$ such that 
a topology $\sf F$ on
$\sf M$ is compatible with $G"$, that is, $\mu \mapsto \mu _h$
is the homomorphism
of $({\sf M,F})$ into itself for each $h \in G"$. Let $\sf F$ be the
topology of convergence for each $E \in Bf(G)$.
Suppose also that $G$ and $G"$ are real Banach manifolds such that
the tangent space $T_eG"$ is dense in $T_eG$, then $TG$ and $TG"$
are also Banach manifolds. Let $\Xi (G")$ denotes
the set of all differentiable vector fields $X$ on $G"$, that is, 
$X$ are sections of the tangent bundle $TG"$. We say that the measure
$\mu $ is continuously differentiable if there exists its tangent mapping
$T_{\phi }\mu _{\phi }(E)(X_{\phi })$ corresponding to the strong
differentiability relative to the Banach structures of the manifolds 
$G"$ and $TG"$. Its differential we denote $D_{\phi }\mu _{\phi }(E)$, 
so $D_{\phi }\mu _{\phi }(E)(X_{\phi })$ is the $\sigma $-additive
real measure by subsets $E\in Af(G, \mu )$ for each $\phi \in G"$ 
and $X\in \Xi (G")$ such that $D\mu (E): TG"\to \bf R$ is continuous
for each $E\in Af(G, \mu )$, 
where $D_{\phi }\mu _{\phi }(E)=pr_2\circ (T\mu )_{\phi }(E)$,
$pr_2: p\times {\bf F}\to \bf F$ is the projection in $TN$, $p\in N$,
$T_pN=\bf F$, $N$ is another real Banach differentiable manifold modelled on 
a Banach space $\bf F$, for a differentiable mapping $V: G"\to N$ 
by $TV: TG"\to TN$ is denoted the corresponding tangent mapping,
$(T\mu )_{\phi }(E):=T_{\phi }\mu _{\phi }(E)$. 
Then by induction $\mu $ is called $n$ times 
continuously differentiable if $T^{n-1}\mu $ is continuously
differentiable such that 
$T^n\mu :=T(T^{n-1}\mu )$, $(D^n\mu )_{\phi }(E)(X_{1, \phi },
...,X_{n, \phi })$ are the $\sigma $-additive real 
measures by $E\in Af(G, \mu )$
for each $X_1,$...,$X_n\in \Xi (G")$, where $(X_j)_{\phi }=:X_{j, \phi }$
for each $j=1,...,n$ and $\phi \in G"$, $D^n\mu : Af(G, \mu )\otimes
\Xi (G")^n\to \bf R$.
\par {\bf 2.1.2.1.} 
~ A canonical closed subset $Q$ of
$X=\bf R^n$ or of the standard separable
Hilbert space $X=l_2({\bf R})$ over $\bf R$
is called a quadrant if it can be given by $Q:=\{ x\in X:
q_j(x)\ge 0 \} $, where $(q_j: j\in \Lambda _Q)$ 
are linearly independent
elements of the topologically adjoint space $X^*$.  Here 
$\Lambda _Q\subset \bf N$ 
(with $card (\Lambda _Q)=k\le n$ when $X=\bf R^n$)
and $k$ is called the index of $Q$.  If $x\in Q$ and exactly $j$ of the $q_i$'s
satisfy $q_i(x)=0$ then $x$ is called a corner of index $j$.  
Since the unitary space $X=\bf C^n$ or the separable Hilbert space
$l_2({\bf C})$ over $\bf C$ as considered over the field $\bf R$
is isomorphic with $X_{\bf R}:=\bf R^{2n}$ or $l_2({\bf R})$ respectively, 
then the above definition also
describes quadrants in $\bf C^n$ and $l_2({\bf C})$ in such sense
(see also \cite{michor}). In the latter case we also consider 
generalized quadrants
as canonical closed subsets which can be given by
$Q:=\{ x\in X_{\bf R}:$ $q_j(x+a_j)\ge 0, a_j\in X_{\bf R},
j\in \Lambda _Q \} ,$ where $\Lambda _Q\subset \bf N$
($card(\Lambda _Q)=k\in \bf N$ when $dim_{\bf R}X_{\bf R}<\infty $).
\par {\bf 2.1.2.2.} If for each open subset $U\subset Q\subset X$ 
a function $f:  Q\to Y$ for Banach spaces $X$ and $Y$ over $\bf R$
has continuous Frech\'et differentials $D^{\alpha }f|_U$ on $U$ 
with $\sup_{x\in U} \| D^{\alpha }f(x) \|_{L(X^{\alpha },Y)} <\infty $ 
for each $0\le \alpha
\le r$ for an integer $0\le r$ or $r=\infty $, 
then $f$ belongs to the class of
smoothness $C^r(Q,Y)$, where $0\le r\le \infty $,
$L(X^k,Y)$ denotes the Banach space of bounded $k$-linear
operators from $X$ into $Y$.
\par {\bf 2.1.2.3.} A differentiable mapping $f:  U\to U'$ is called a 
diffeomorphism if 
\par $(i)$ $f$ is bijective and there exist continuous
$f'$ and $(f^{-1})'$, where $U$ and $U'$ 
are interiors of quadrants $Q$
and $Q'$ in $X$. 
\par In the complex case we consider bounded generalized 
quadrants $Q$ and $Q'$ in $\bf C^n$ or $l_2({\bf C})$
such that they are domains with piecewise 
$C^{\infty }$-boundaries and we impose additional conditions on
the diffeomorphism $f$:
\par $(ii)$ ${\bar \partial }f=0$ on $U$, 
\par $(iii)$ $f$ and all its strong (Frech\'et) differentials (as multilinear
operators) are bounded on $U$, 
where $\partial f$ and ${\bar \partial }f$
are differential $(1,0)$ and $(0,1)$ forms respectively,
$d=\partial +{\bar \partial }$ is an exterior derivative.
In particular for 
$z=(z^1,...,z^n)\in \bf C^n$, $z^j\in \bf C$, $z^j=x^{2j-1}+ix^{2j}$
and $x^{2j-1}, x^{2j}\in \bf R$ for each $j=1,...,n,$
$i=(-1)^{1/2}$, there are expressions:
$\partial f:=\sum_{j=1}^n(\partial f/\partial z^j)dz^j$,
${\bar \partial }f:=\sum_{j=1}^n(\partial f/\partial {\bar z}^j)d{\bar z}^j$.
In the infinite-dimensional case there are equations:
$(\partial f)(e_j)=\partial f/\partial z^j$
and $({\bar \partial }f)(e_j)=\partial f/\partial {\bar z}^j$, 
where $\{ e_j: {j\in \bf N} \} $ is the standard orthonormal base
in $l_2({\bf C})$, $\partial f/\partial z^j=(\partial f/\partial x^{2j-1}
-i\partial f/\partial x^{2j})/2$,
$\partial f/\partial {\bar z}^j=(\partial f/\partial x^{2j-1}
+i\partial f/\partial x^{2j})/2$.
\par Cauchy-Riemann Condition $(ii)$ means that $f$ on $U$ is the holomorphic
mapping.
\par {\bf 2.1.2.4.}  A complex manifold $M$ with corners is
defined in the usual way:  it is a metric separable space
modelled on $X=\bf C^n$ or $X=l_2({\bf C})$ 
and is supposed to be of class $C^{\infty }$.  Charts on $M$ are
denoted $(U_l, u_l, Q_l)$, that is $u_l:  U_l\to u_l(U_l) \subset Q_l$ are
$C^{\infty }$-diffeomorphisms, $U_l$ are open in $M$, $u_l\circ
{u_j}^{-1}$ are biholomorphic from
domains $u_j(U_l\cap U_j)\ne
\emptyset $ onto $u_l(U_l\cap U_j)$ 
(that is $u_j\circ u_l^{-1}$ and 
$u_l\circ u_j^{-1}$ are holomorphic and bijective)
and $u_l\circ u_j^{-1}$ 
satisfy conditions $(i-iii)$ from \S 2.1.2.3, $\bigcup_jU_j=M$
(see also Note 2.2.3).
\par A point $x\in M$ is called a corner of index $j$
if there exists a chart $(U,u,Q)$ of $M$ with $x\in U$ and $u(x)$ is of index
$ind_M(x)=j$ in $u(U)\subset Q$.  The set of all corners of index $j\ge 1$ is
called the border $\partial M$ of $M$, $x$ is called an inner point of $M$ if
$ind_M(x)=0$, so $\partial M=\bigcup_{j\ge 1}\partial ^jM$, where
$\partial ^jM:=\{ x\in M:  ind_M(x)=j \} $.  
\par For the real manifold with corners on the connecting mappings
$u_l\circ u_j^{-1}\in C^{\infty }$ 
of real charts is imposed only Condition $2.1.2.3(i)$.
\par {\bf 2.1.2.5.} A subset $Y\subset M$ is called a 
submanifold with corners of $M$ if for each
$y\in Y$ there exists a chart $(U,u,Q)$ of $M$ 
centered at $y$ (that is $u(y)=0$
) and there exists a quadrant $Q'\subset {\bf C^k}$ or in $l_2({\bf C})$
such that
$Q'\subset Q$ and $u(Y\cap U)=u(U)\cap Q'$.  A submanifold with corners $Y$ of
$M$ is called neat, if the index in $Y$ of each $y\in Y$ coincides with its
index in $M$. 
\par  Analogously for real manifolds with corners for $\bf R^k$ and
$\bf R^n$ or $l_2({\bf R})$ instead of $\bf C^k$ and $\bf C^n$
or $l_2({\bf C})$.
\par {\bf 2.1.2.6.} Henceforth, the term a complex manifold
$N$ modelled on $X=\bf C^n$ or $X=l_2({\bf C})$ means a metric separable
space supplied with an atlas $\{ (U_j,\phi _j): j\in \Lambda _N \} $
such that:
\par $(i)$ $U_j$ is an open subset of $N$ for each $j\in \Lambda _N$
and $\bigcup_{j\in \Lambda _N}U_j=N$, where $\Lambda _N\subset \bf N$;
\par $(ii)$ $\phi _j: U_j\to \phi _j(U_j)\subset X$ are
$C^{\infty }$-diffeomorphisms, where $\phi _j(U_j)$ are 
$C^{\infty }$-domains in $X$;
\par $(iii)$ $\phi _j\circ \phi _m^{-1}$ is a biholomorphic mapping
from $\phi _m(U_m\cap U_j)$ onto
$\phi _j(U_m\cap U_j)$ while $U_m\cap U_j\ne \emptyset $. 
When $X=l_2({\bf C})$ it is supposed, that
$\phi _j\circ \phi _m^{-1}$ are Frech\'et (strongly)
$C^{\infty }$-differentiable.
\par {\bf 2.1.3.1.}  Let $X$ be either the standard separable
Hilbert space $l_2=l_2({\bf C})$ over the field $\bf C$ of complex numbers
or $X=\bf C^n$.  Let $t\in \bf
N_o$ $:={\bf N}\cup \{ 0\}$, ${\bf N}:=\{ 1,2,3,... \}$ 
and $W$ be a domain with a continuous piecewise
$C^{\infty }$-boundary $\partial W$
in $\bf R^{2m}$,
$m\in \bf N$, that is $W$ is a $C^{\infty }$-manifold with corners
and it is a canonical closed
subset of $\bf C^m$, $cl(Int(W))=W$, where $cl(V)$ denotes the closure of $V$,
$Int(V)$ denotes the interior of $V$ in the corresponding 
topological space.  As
usually $H^t(W,X)$ denotes the Sobolev space of functions 
$f:  W\to X$ for which
there exists a finite norm
\par $\| f\|_{H^t(W,X)}:=(\sum_{|\alpha |\le t}{\|
D^{\alpha }f\|^2}_{L^2 (W,X)})^{1/2}<\infty $, \\
where $f(x)=(f^j(x):  j\in {\bf
N})$, $f(x)\in l_2$, $f^j(x)\in \bf C$, $x\in W$, 
\par ${\| f\|^2}_{L^2
(W,X)}:=\int_W {\| f(x)\|^2}_X\lambda (dx)$, $\lambda $ is the Lebesgue measure
on $\bf R^{2m}$, $\| z\|_{l_2}:= (\sum_{j=1}^{\infty }|z^j|^2)^{1/2}$, $z=(z^j:
j\in {\bf N})\in l_2$, $z^j\in \bf C$. Then $H^{\infty }(W,X):=
\bigcap_{t\in \bf N}H^t(W,X)$ is the uniform space with the uniformity
given by the family of norms $\{ \| f \|_{H^t(W,X)}: t\in {\bf N} \}$.
\par {\bf 2.1.3.2.}  Let now $M$ be a
compact Riemann or complex
$C^{\infty }$-manifold with corners with 
a finite atlas $At(M):=\{ (U_i, \phi _i,
Q_i); i\in \Lambda _M \} $, where $U_i$ are open in $M$, 
$\phi _i:  U_i\to \phi
_i(U_i)\subset Q_i\subset \bf R^m$ (or it is a subset in $\bf C^m$)
are diffeomorphisms
(in addition holomorphic respectively as in \S 2.1.2.3), $(U_i, \phi _i)$ are
charts, $i\in \Lambda _M \subset \bf N$.  
\par Let also $N$ be a
separable complex metrizable manifold with corners 
modelled either on $X=\bf C^n$
or on $X=l_2({\bf C})$ respectively. Let 
$(V_i, \psi _i, S_i)$ be
charts of an atlas $At(N):=\{ (V_i, \psi _i, S_i):  i\in \Lambda _N \} $
such that
$\Lambda _N\subset \bf N$ and $\psi _i:  V_i\to \psi _i(V_i)\subset 
S_i\subset X$ are
diffeomorphisms, $V_i$ are open in $N$, $\bigcup_{i\in \Lambda _N}V_i=N$.
We denote by $H^t(M,N)$ the Sobolev
space of functions $f:  M\to N$ for which $f_{i,j}\in H^t(W_{i,j},X)$ for each
$j\in \Lambda _M$ and $i\in \Lambda _N$ for a domain $W_{i,j}\ne \emptyset $ of
$f_{i,j}$, where $f_{i,j}:=\psi _i\circ f\circ {\phi _j}^{-1}$, and
$W_{i,j}=\phi
_j(U_j\cap f^{-1}(V_i))$ are canonical closed subsets of $\bf R^m$
(or $\bf C^m$ respectively).  The
uniformity in $H^t(M,N)$ is given by the following base $\{ (f,g)\in
(H^t(M,N))^2:  \sum_{i\in \Lambda _N, j\in \Lambda _M} {\| f_{i,j}-
g_{i,j}\|^2}_{H^t(W_{i,j}, X)}<\epsilon \}$, where $\epsilon >0$, $W_{i,j}$
are domains of $(f_{i,j}- g_{i,j}).$ For $t=\infty $ as usually
$H^{\infty }(M,N):=\bigcap_{t\in \bf N} H^t(M,N)$.
\par {\bf 2.1.3.3.}  For two complex manifolds $M$ and $N$ 
with corners let 
${\sf O}_{\Upsilon }(M,N)$ denotes a space
of continuous mappings $f: M\to N$ such that
for each $f$ there exists a partition $Z_f$ of $M$ with the help of
a real $C^{\infty }$-submanifold ${M'}_f$, 
which may be with corners, such that its codimension 
over $\bf R$ in $M$ is $codim_{\bf R}{M'}_f=1$
and $M\setminus {M'}_f$ is a disjoint union of open complex submanifolds
$M_{j,f}$ possibly with corners 
with $j=1,2,...$ such that each restriction $f|_{M_{j,f}}$
is holomorphic with all its derivatives bounded on $M_{j,f}$.
For a given partition $Z$ (instead of $Z_f$) and the corresponding
$M'$ the latter subspace of continuous piecewise holomorphic mappings
$f: M\to N$ is denoted by ${\sf O}(M,N;Z)$.
The family $\{ Z \} $ of all such partitions is denoted
$\Upsilon $.
That is ${\sf O}_{\Upsilon } (M,N)=str-ind_{\Upsilon }(M,N;Z)$.
Let also ${\sf O}(M,N)$ denotes the space
of holomorphic mappings $f: M\to N$,
$Diff^{\infty }(M)$ denotes a group of 
$C^{\infty }$-diffeomorphisms of $M$ and $Diff^{\sf O}_{\Upsilon }(M):=
Hom(M)\cap {\sf O}_{\Upsilon }(M,M)$, where $Hom(M)$ is a group
of homeomorphisms.
\par Let $A$ and $B$ be two complex manifolds with corners
such that $B$ is a submanifold of $A$. Then $B$ is called
a strong $C^r([0,1]\times A,A)$-retract 
(or $C^r([0,1],{\sf O}_{\Upsilon }(A,A))$-retract) of  $A$ if there exists
a mapping $F: [0,1]\times A\to A$ such that $F(0,z)=z$ for each
$z\in A$ and $F(1,A)=B$ and $F(x,A)\supset B$ for each $x\in
[0,1]:=\{ y: 0\le y \le 1, y\in {\bf R} \} $,
$F(x,z)=z$ for each $z\in B$ and $x\in [0,1]$,
where $F\in C^r([0,1]\times A,A)$ or $F\in C^r([0,1],{\sf O}_{\Upsilon }
(A,A))$ respectively, $r\in [0,\infty )$, $F=F(x,z)$, $x\in [0,1]$,
$z\in A$. Such $F$ is called the retraction. In the case of $B=\{ a_0 \} $,
$a_0\in A$ we say that $A$ is $C^r([0,1]\times A,A)$-contractible
(or $C^r([0,1], {\sf O}_{\Upsilon }(A,A))$-contractible correspondingly).
Two maps $f: A\to E$ and $h: A\to E$ are called 
$C^r([0,1]\times A,E)$-homotopic (or 
$C^r([0,1],{\sf O}_{\Upsilon }(A,E))$-homotopic ) if there
exists $F\in C^r([0,1]\times A,E)$ (or $F\in C^r([0,1],{\sf O}_{\Upsilon }
(A,E))$ respectively) such that $F(0,z)=f(z)$ and $F(1,z)=h(z)$
for each $z\in A$, where $E$ is also a complex manifold.
Such $F$ is called the homotopy.
\par Let $M$ be a complex manifold with corners
satisfying the following conditions:
\par $(i)$ it is compact;
\par $(ii)$ $M$ is a union of two closed complex submanifolds $A_1$ and $A_2$
with corners,
which are canonical closed subsets in $M$
with $A_1\cap A_2=\partial A_1\cap \partial A_2=:A_3$ 
and a codimension over $\bf R$ of $A_3$ in $M$ is $codim_{\bf R}A_3=1$;
\par $(iii)$ a marked point $s_0$ is in $A_3$;
\par $(iv)$ $A_1$  and $A_2$ are 
$C^0([0,1],{\sf O}_{\Upsilon }(A_j,A_j))$-contractible into a marked point
$s_0\in A_3$ by mappings $F_j(x,z),$
where either $j=1$ or $j=2$.
\par  We consider all finite partitions $Z:=\{ M_k:  k\in \Xi _Z\}$ of $M$
such that $M_k$ are complex submanifolds (of $M$), which may be
with corners and
$\bigcup_{k=1}^sM_k=M$, $\Xi _Z=\{ 1,2,...,s \}$, $s\in \bf N$ 
depends on $Z$, $M_k$ are
canonical closed subsets of $M$.  We denote by $\tilde diam(Z):=\sup_k(diam
(M_k))$ the diameter of the partition $Z$, where $diam (A)=\sup_{x, y\in A}
|x-y|_{\bf C^n}$ is a diameter of a subset $A$ in $\bf C^n$, since
each finite-dimensional manifold $M$ can be embedded into $\bf C^n$
with the corresponding $n\in \bf N$.  
We suppose
also that $M_i\cap M_j\subset M'$ and $\partial M_j\subset M'$ for each $i\ne
j$, where $M'$ is a closed $C^{\infty }$-submanifold (which may be
with corners) in $M$ with
the codimension $codim_{\bf R}(M')=1$ of $M'$ in $M$, $M'=\bigcup_{j\in \Gamma
_Z}{M'}_j$, ${M'}_j$ are $C^{\infty }$-submanifolds of $M$, 
$\Gamma _Z$ is a finite
subset of $\bf N$. 
\par  We denote by $H^t(M,N;Z)$ a space of continuous functions
$f:  M\to N$ such that $f|_{(M\setminus M')}\in H^t(M\setminus M',N)$ and
$f|_{[Int(M_i)\cup (M_i\cap {M'}_j)]}\in H^t(Int(M_i)\cup (M_i\cap {M'}_j),N)$,
when $\partial M_i\cap {M'}_j\ne \emptyset $, $h^Z_{Z'}:  H^t(M,N;Z)\to
H^t(M,N;Z')$ are embeddings for each $Z\le Z'$ in $\Upsilon $.
\par The ordering
$Z\le Z'$ means that each submanifold $M_i^{Z'}$
from a partition $Z'$ either belongs to the family
$(M_j:  j=1,...,k)=(M_j^Z: j=1,...,k)$ for $Z$ or 
there exists $j$ such that $M_i^{Z'}\subset M_j^Z$ and
$M_j^Z$ is a finite union of $M_l^{Z'}$ for which $M_l^{Z'}\subset
M_j^Z$. Moreover, these $M_l^{Z'}$ are submanifolds (may be
with corners) in $M_j^Z$. 
\par Then we consider the following uniform space
$H^t_p(M,N)$ that is the strict inductive limit $str-ind \{ H^t(M,N;Z);
h^{Z'}_{Z}; \Upsilon \} $ (the index $p$ reminds about the procedure of
partitions), where $\Upsilon $ is the directed family of all such $Z$, 
for which $\lim_{\Upsilon }\tilde diam(Z)=0$.  
\par {\bf 2.1.4.}  Let now $s_0$ be the marked
point in $M$ such that $s_0\in A_3$ (see \S 2.1.3.3) 
and $y_0$ be a marked point in the manifold $N$.
\par $(i).$ Suppose that $M$ and $N$ are connected.
\par Let
$H^t_p(M,s_0;N,y_0):=\{ f\in H^t(M,N)| f(s_0)=y_0 \} $ denotes the closed
subspace of $H^t(M,N)$ and $\omega _0$ be its element such that $\omega
_0(M)= \{ y_0 \} $, where $\infty \ge t\ge m+1$, $2m=dim_{\bf R}M$
such that $H^t\subset C^0$ due to the Sobolev embedding theorem.  
The following subspace $ \{ f: f\in H^{\infty }_p(M,s_0;N,y_0),\mbox{ }
\mbox{ } {\bar \partial }f=0 \} $ is isomorphic with
${\sf O}_{\Upsilon }(M,s_0;N,y_0)$, since $f|_{(M\setminus M')}
\in H^{\infty }(M\setminus M',N)=
C^{\infty }(M\setminus M',N)$ and ${\bar \partial }f=0$.
\par Let as usually $A\vee B:=A\times \{
b_0\}\cup \{ a_0\}\times B\subset A\times B$ 
be the wedge sum of 
pointed spaces $(A,a_0)$ and $(B,b_0)$, where $A$
and $B$ are topological spaces with marked points $a_0\in A$ and $b_0\in B$.
Then the wedge combination $g\vee f$ 
of two elements $f, g\in H^t_p(M,s_0;N,y_0)$ is
defined on the domain $M\vee M$ (see also
Chapters 0-2 \cite{swit} and Example 2.4.12 \cite{eng}).
\par The spaces ${\sf O}_{\Upsilon }(J,A_3;N,y_0):=\{ f\in {\sf O}_{\Upsilon }
(J,N): f(A_3)=\{ y_0 \} \} $ 
have the manifold structure and have embeddings into
${\sf O}_{\Upsilon }(M,s_0;N,y_0)$ due to Condition 
$2.1.3.3(ii)$, where either
$J=A_1$ or $J=A_2$.
This induces the following embedding
$\chi ^*: {\sf O}_{\Upsilon }(M\vee M,s_0\times s_0;N,y_0)\hookrightarrow
{\sf O}_{\Upsilon }(M,s_0;N,y_0)$.
Therefore $g\circ f:=\chi ^*(f\vee g)$ is the composition in
${\sf O}_{\Upsilon }(M,s_0;N,y_0)$.  
\par The space $C^{\infty }(M,N)$ is dense in $C^0(M,N)$
and there is the inclusion ${\sf O}_{\Upsilon }(M,N)\subset
H^{\infty }_p(M,N)$. Let $M_{\bf R}$ be the Riemann manifold generated by
$M$ considered over $\bf R$. Then $Diff^{\infty }_{s_0}(M_{\bf R})$
is a group of $C^{\infty }$-diffeomorphisms $\eta $ of $M_{\bf R}$ 
preserving the marked point $s_0$, that is $\eta (s_0)=s_0$.
There exists the following equivalence relation
$R_{\sf O}$ in ${\sf O}_{\Upsilon }(M,s_0;N,y_0)$:
$fR_{\sf O}h$ if and only if there exist nets
$\eta _n\in Diff^{\infty }_{s_0}(M_{\bf R})$, 
also $f_n$ and $h_n\in H^{\infty }_p(M,s_0;N,y_0)$ with $\lim_{n}f_n=f$ 
and $\lim_{n}h_n=h$ such that $f_n(x)=h_n(\eta _n(x))$ for each
$x\in M$ and $n\in \omega $, where $\omega $ is an ordinal,
$f, h \in {\sf O}_{\Upsilon }(M,s_0;N,y_0)$ and converegence is considered
in $H^{\infty }_p(M,s_0;N,y_0)$. 
\par The quotient space ${\sf O}_{\Upsilon }(M,s_0;N,y_0)/R_{\sf O}=:
(S^MN)_{\sf O}$ is called the loop semigroup. It will be shown later, that
$(S^MN)_{\sf O}$ has a structure of topological Abelian monoid
with the cancellation property.
Applying the A. Grothendieck procedure \cite{langal,swan}
to $(S^MN)_{\sf O}$ we get a loop group $(L^MN)_{\sf O}$.
For the spaces $H^t_p(M,s_0;N,y_0)$ the corresponding equivalence 
relations are denoted $R_{t,H}$, the group semigroups are denoted 
by $(S^M_{\bf R}N)_{t,H}$, the loop groups are denoted
by $(L^M_{\bf R}N)_{t,H}$.
\par {\bf 2.2. Propositions. (1).} {\it Let $f$ be a diffeomorphism
satisfying Conditions $2.1.2.3(i-iii),$ then there are
neighbourhoods $V$ and $V'$ of bounded generalized
quadrants $Q$ and $Q'$ such that
$f$ has an extension $f$ to the holomorphic diffeomorphism 
of $V$ with $V'$.}
\par {\bf (2).} {\it If $f: U\to \bf C$ is a mapping satisfying Conditions
$2.1.2.3(ii,iii),$ then there exists a neighbourhood $V$ of $Q$
such that $f$ has a holomorphic extension onto $V$.}
\par {\bf Proof. (1). } From \cite{seel,touger} it follows, that
$f$ has a $C^{\infty }$-extension $h$ onto an open domain $W$
with a $C^{\infty }$-boundary such that $W\supset Q$, since
due to Condition $2.1.2.3 (iii)$
each partial derivative \\
$D^{\alpha }f:= \partial ^{|\alpha |}f/\partial (z^1)^{\alpha _1}
...
\partial (z^n)^{\alpha _n}\partial ({\bar z}^1)^{\alpha _{n+1}}...
\partial ({\bar z}^n)^{\alpha _{2n}}$ has a bounded continuous extension
onto $Q$ due to the line integration of $D^{\beta }f$ 
along $C^{\infty }$-curves, where $\beta =\alpha +e_j$, $\alpha =
(\alpha _1,...,\alpha _{2n})$, $0\le \alpha _j\in \bf Z$,
$j=1,...,2n$, since $Q$ is simply connected. 
Therefore, ${\bar \partial }f=0$ on $Q$.
\par If $F'(z)=f(z)$ and $f(z)$ does not satisfy the Cauchy-Riemann
conditions for $z\in \partial Q$, then $F(z)$ also does not satisfy
them.
In view of Corollary 3.2.3, Conjecture in Exer. 3.2 
and Exer. 1.28 \cite{henle} (see also references therein)
there exists a 
piecewise holomorphic function $\Phi $ on a bounded neighbourhood $V$ of $Q$
in $\bf C^n$ such that on
$\partial Q$ it satisfies the jump condition $\Phi ^-(z)=\Phi ^+ (z)+
g(z)$ for each $z\in \partial Q$, $\Phi |_{V\setminus Q}$ 
and $\Phi |_{Int(Q)}$
are holomorphic and bounded together with each partial derivative,
where $(\alpha )$ $g\in H^{\infty }_p(\partial Q)$,
$(\beta )$ $g|_S\in C^{\infty }(S)$ for each $C^{\infty }$-submanifold $S$
(without corners) in $\partial Q$,
$(\gamma )$ $\int_{\partial Q}g\wedge {\bar \partial }p=0$ for each
$C^{\infty }_{(n,n-2)}$-form $p$ in a neighbourhood of $\partial Q$ 
in $\bf C^n$,
restrictions $\Phi ^+(z)|_{\partial Q}$ are taken as limits 
$\lim_{a\to z, a\in U}\Phi (a)$ and restrictions
$\Phi ^-(z)|_{\partial Q}$ are taken as 
$\lim_{a\to z, a\in V\setminus Q}\Phi (a)$,
which supposed to be existent, since
this is based on the existence of a solution $u$ of
the ${\bar \partial } $-equation ${\bar \partial }u=f$
and the complex conjugate of $\partial \bar u$
is equal to ${\bar \partial }u$ for the corresponding differential forms,
where $V\subset W$.
If $g$ satisfies Conditions $(\alpha ,\beta )$ only, 
then $\Phi |_{V\setminus Q}$
and $\Phi |_{Int(Q)}$ are of class $C^{\infty }$
and bounded together with each partial derivative.
Indeed, there exists an increasing sequence $W_j$ of 
$C^{\infty }$-subregions in $Q$ such that
$(i)$ $\sup_k diam(E_{ k,j})<j^ {-2}$;
$(ii)$ $\eta _j: \partial Q_j\to \partial W_j$ are 
${\sf O}_{\Upsilon }$-diffeomorphisms (that is homeomorphisms of class
${\sf O}_{\Upsilon }$) with $\eta _j|_{\partial Q_j\cap \partial W_j}=id$;
$(iii)$ $g_j$ are $C^{\infty }$-functions satisfying Condition $(\gamma )$
on $\partial W_j$ (as $g$ on $\partial Q$) and
$\lim_{j\to \infty }g_j\circ \eta _j=g$;
$(iv)$ $E_{k,j}$ are connected components of 
$\partial W_j\setminus \partial Q_j$ such that $E_{k,j}\subset
(\partial ^2Q)^{1/j}$, 
where $A^{\epsilon }:= \{ y\in {\bf C^n}:$
$inf_{z\in A} |y-z|<\epsilon \} $ is an $\epsilon $-enlargement
of a subset $A$ in $\bf C^n$, $W_j\subset W_{j+1}$ for each $j$. 
In view of Monthel Theorem 2.14.1 \cite{heins}
there exists a sequence of functions
$\{ \Phi _j:$ $j=l+1, l+2,... \} $ 
satisfying the same conditions for $(W_j,g_j)$
as $\Phi $ for $(Q,g)$ such that $\Phi _j$ converges uniformly
on $\bigcup_{j=l+1}^{\infty }W_j$ and on $V\setminus Q$, where $V$ 
is a bounded neighbourhood of $Q$ in $\bf C^n$. Since $\bigcup_{j=l+1}^{\infty }
W_j\supset (Q\setminus {\tilde \partial }^2Q)$, 
$Q\cap (\bigcap _{j=l+1}^{\infty }(V\setminus Int (W_j))=\partial Q$
and due to the solution of the $\bar \partial $-equation on $Q$ and on
$V\setminus Int (Q)$ there exists the sequence $\{ \Phi _j: j>l \} $
which converges to the desired $\Phi $ on $V$, 
where ${\tilde \partial }^pQ:=\bigcup_{j, j\ge p}\partial ^jQ$,
$p\in \bf N$ (see \S 2.1.2.4).
In the particular case of a $C^1$-domain $Q$ in $\bf C$ see also
\cite{msh1,msh2}.
The line integration of $\Phi $ along $C^{\infty }$-curves
produces a continuous $\Phi $ with the jump condition for $\partial \Phi $
instead of $\Phi $ due to Cauchy Integral Theorem
in the Stronger Form 6.4.1 \cite{heins} and \cite{henle}, 
since $Q$ is simply connected (see also \S 3.5 \cite{moret}).
\par Therefore, there exists
$C^{\infty }$-function $v$ on $V$ such that $v|_{Q}=0$ and 
${\bar \partial }v= -{\bar \partial }h$ on $V$. That to construct such
$v$ let $v=v_1-v_2-v_3$, where $v_1|_{Q}=(v_2+v_3)|_{Q}$, 
$v_1$ is a $C^{\infty }$-function on $U=Q\setminus \partial Q$ 
and on $V\setminus Q$,
$v_2$ and $v_3$ are holomorphic on
$U$ and $V\setminus Q$, $v_1^+|_{\partial Q}=f^+|_{\partial Q}$, 
$v_2^+|_{\partial Q}=0$, $\partial v_2^-|_{\partial Q}=
\partial f^+|_{\partial Q}$,
$v_3^+|_{\partial Q}=f^+|_{\partial Q}$, 
$v_3^-|_{\partial Q}=0$.
In view of the theorem about maximum principle of the modulus of
a holomorphic function on a simply connected domain the condition
$v_2^{+}|_{\partial Q}=0$ implies $v_2|_Q=0$ (see \S 3.6 \cite{moret}).
Therefore, $\bar \partial v|_Q=\bar \partial
v_1|_Q$ and $\bar \partial (h+v)|_V=0$ and $(h+v)|_Q=f$.
Then $u:=h+v$ is a holomorphic extension of $f$ on $V$.
When $n\ge 2$ Hartogs Theorem 1.2.2 \cite{henle} can be 
used instead of Conjecture in Exer. 3.2, that is for the constructed above 
$u$ holomorphic on $V\setminus {\tilde \partial }^2Q$ there exists
a function $v$ holomorphic on $V$ such that 
$v|_{V\setminus {\tilde \partial }^2Q}=u|_{V\setminus {\tilde \partial }^2Q}$,
since $V\setminus {\tilde \partial }^2Q$ is connected.
For $Q\subset {\bf C^1}=X$ there exists a strictly pseudoconvex open set 
$D\subset \bf C^2$ such that $X\cap {\bar D}=Q$, hence by 
Theorem 3.6.8 \cite{henle} there exists $F\in {\sf O}(D)\cap C^0({\bar D})$ 
such that $F|_Q=f$, consequently, the case of $n=1$ can be reduced to the 
case $n=2$ with a subsequent restriction of a resulting function on $X\cap V$.
\par Piecewise $\partial Q$ can be written in local coordinates
in a neighbourhood $W_z$  of $z\in \partial Q$
as $x^1=0$, where $z=(z^1,z^2,...)$, $z^j=x^{2j-1}+ix^{2j}$.
Therefore, $f$ has also such extension for $Q\subset l_2({\bf C})$.
\par On the other hand, the Jacobian $J_f$ is not equal to zero 
everywhere on the closed set $Q$ in $X=\bf C^n$
or $X=l_2({\bf C})$. 
In view of the strong $C^{\infty }$-differentiability of $f$
there exists
an open set $U$ in $X$ such that $Q\subset U\subset V$
and $J_f\ne 0$ everywhere on $U$. Then we take $U'=f(U)$,
hence $Q'\subset U'$ and $U'$ is open in $X$.
\par {\bf (2).} The proof of this case is analogous to that of
Section $(1)$ omitting the last paragraph from it.
\par {\bf 2.2.3. Note.} From the proof it follows 
that these propsitions are true, 
when Condition $2.1.2.3(iii)$ is fulfilled up to
the second differentials on $U$, 
since if $f\in C^1(V,{\bf C})$ and $f$ satisfies
Cauchy-Riemann conditions on $V$ it follows that $f\in {\sf O}(V,{\bf C})$.
These propositions also justify the notion of $f$ to be holomorphic on
$Q$ as a restriction $f|_Q$ of $f$ holomorphic on $V$.
\par {\bf 2.3.1. Lemma.} {\it If $M$ is a complex manifold 
modelled on $X=\bf C^n$ or $X=l_2({\bf C})$ 
with an atlas $At(M)=\{ (V_j,\phi _j): j \}$,
then there exists an atlas $At'(M)= \{ (U_k,u_k,Q_k): k \} $
which refines $At(M)$, where $(V_j,\phi _j)$ are usual charts
with diffeomorphisms $\phi _j: V_j\to \phi _j(V_j)$ such that 
$\phi _j(V_j)$ are $C^{\infty }$-domains in $\bf C^n$ and 
$(U_k,u_k,Q_k)$ are charts corresponding to quadrants $Q_k$ in $\bf C^n$
or $l_2({\bf C})$.}
\par {\bf Proof.} The covering $\{ V_j: j \} $ of $M$ has a refinement
$\{ W_l: l \} $ such that for each $j$ there exists $l=l(j)$
with $W_l\subset V_j$ so that each $\phi _j(W_l)$ is a simply connected 
region in $\bf C^n$ or $l_2({\bf C})$ which is not the whole space.
We choose $W_l$ such that 
\par $(i)$ either $W_l\cap \partial M=\emptyset $
or $W_l\cap \partial M$ is open in $\partial M$;
\par $(ii)$ $ \{ \pi _k(z)=z^k:$ $z\in \phi _j(W_l) \} =:E_{l,k}$,
$z\in X$, $\pi _k:$ $X\to \bf C$ are canonical projections
associated with the standard orthonormal base $\{ e_j:$ $j \} $ in  $X$,
$E_{l,k}$ are $C^{\infty }$-regions in $\bf C$, $\phi _j(W_l)=
\prod_{k=1}^nE_{l,k}$ for $X=\bf C^n$, or $\pi _J(\phi _j(W_l))=
\prod_{k\in J}E_{l,k}$ for each finite subset $J$ in $\bf N$
and the corresponding projection $\pi _J: X\to sp_{\bf C} \{ e_k:$ $k\in J \} $.
In view of the Riemann Mapping Theorem 
for each $E_{l,k}$ there exists holomorphic 
diffeomorphism $q_{l,k}$ either onto $B^{-}_r:=\{ z\in {\bf C} :$ $|z|<r \} $
or $F_r:=\{ z\in {\bf C}:$ $|z|<r, x^1\ge 0 \} $
(see \S 2.12 in \cite{heins}). The latter case appears
while treatment of $\pi _k(\phi _j(W_l\cap \partial M))\ne \emptyset $
(see \S 10.5.2 \cite{sanger} and \S 12 \cite{heins}). Slightly shrinking
covering if necessary we can choose $\{ W_l : l \} $ such
that  each $q_{l,k}$ and its derivatives are bounded on $E_{l,k}$.
In view of Central Theorem from \S 6.3 \cite{heins}
$q_{l,k}$ are boundary preserving maps. In view of Chapter 13 \cite{moret}
$B^{-}_r$ and $F_r$ have finite atlases with charts corresponding to quadrants
(see \S 2.1.2 and \S 2.2).
\par {\bf 2.3.2} {\bf Note.} Vice versa there are complex manifolds 
with corners, which are not usual complex manifolds, for example,
canonical closed domains $F$ in $\bf C^n$ with piecewise 
$C^{\infty }$-boundary, which is not of class $C^1$.
Since each complex manifold $G$ has $\partial G$ of class $C^{\infty }$
by Definition 2.1.2.6.
\par {\bf 2.4. Lemma.} {\it ${\sf O}_{\Upsilon }(M,N)$ from \S 2.1.3.3
is an infinite-dimensional complex manifold dense in $C^0(M,N)$,
when $M$ has $dim_{\bf C}M\ge 1$. Moreover, there exists its tangent bundle
$T{\sf O}_{\Upsilon }(M,N)={\sf O}_{\Upsilon }(M,TN).$ If $N=\bf C^n$ or
$N=l_2({\bf C})$, then ${\sf O}_{\Upsilon }(M,N)$ is an infinite-dimensional
topological vector space over $\bf C$.}
\par {\bf Proof.} The connecting mappings $\phi _j\circ \phi _k^{-1}$
of charts $(U_j,\phi _j)$ and $(U_k,\phi _k)$ with $U_j\cap U_k\ne 
\emptyset $ are holomorphic on the corresponding domains
$\phi _k(U_j\cap U_k)$. 
For each submanifold $M_j$ in $M$ we have $T{\sf O}(M_j,N)=
{\sf O}(M_j,TN)$ \cite{eliass}.
For each $f\in {\sf O}_{\Upsilon }(M,N)$ there exists a partition
$Z_f$ of $M$ such that $f|_{M_{j,f}}\in {\sf O}(M_{j,f},N)$ for each 
submanifold $M_{j,f}$ with corners  defined by $Z_f$.
In accordance with
\S 2.1.3.2 and \S 2.1.4 the topology of ${\sf O}_{\Upsilon }(M,N)$
is the compact-open topology, hence $\phi _j\circ \phi _k^{-1}$
induce connecting mappings $(\phi _j^{-1}\circ \phi _k)^*$ of the corresponding
charts in ${\sf O}_{\Upsilon }(M,N)$ such that 
$(\phi _j^{-1}\circ \phi _k)^*(f(z)):=f\circ (\phi _j^{-1}\circ \phi _k)(z)$
for each $z\in U_j\cap U_k$ such that its Frech\'et derivatives are
the following $[\partial (\phi _j^{-1}\circ \phi _k)^*(f)/\partial f].h=
(\phi _j^{-1}\circ \phi _k)^*(h)$ and 
$[\partial (\phi _j^{-1}\circ \phi _k)^*(f)/\partial \bar f].h=0$,
where $h$ are vectors in $T_f{\sf O}_{\Upsilon }(U_j\cap U_k,N)$.
Therefore,
$T{\sf O}_{\Upsilon }(M,N)={\sf O}_{\Upsilon }(M,TN)$, since
${\sf O}_{\Upsilon }(M,N)$ is the complex manifold (certainly of class 
$C^{\infty }$). 
\par In particular ${\sf O}_{\Upsilon }(M,Y)$ is a topological vector space
over $\bf C$ for $Y=\bf C^n$ or $Y=l_2({\bf C})$.
It remains to prove that ${\sf O}_{\Upsilon }(M,N)$ is infinite-dimensional
and dense in $C^0(M,N)$. This follows from 
Corollary 3.2.3, Exer. 1.28 and Conjecture in Exer. 3.2 
\cite{henle}, since for each
quadrant $Q$ and a given function $s$ on $\partial Q$,
which is a restriction $q|_{\partial Q}$ of a holomorphic
function $q$ on a neighbourhood of $\partial Q$ in $\bf C^m$,
$m=dim_{\bf C}M$,
there exists a space of functions $u: W\to \bf C^n$
such that $u|_{Int(Q)}$ and $u|_{W\setminus Q}$ 
are holomorphic and bounded together with each
partial derivative, where $W$ is
an open ball in $\bf C^n$ containing $Q$ and such that
$u^+(z)-u^-(z)=s(z)$ for each $z\in \partial Q$.
Then using Cauchy integration along $C^{\infty }$-curves
we construct a space of continuous functions 
$f: W\to \bf C^n$ holomorphic on $U$ and $W\setminus Q$
with prescribed  $(\partial f)^-(z)-(\partial f)^{+}(z)$ 
for each $z\in \partial Q$ and analogously for $f\in C^l$
with jump conditions for higher order derivatives.
In the case of $n>1$ there also can be used holomorphic extension
of holomorphic functions from proper complex submanifolds $K$
of $\partial Q$ (see Theorems 2(b) and 3(b) in \cite{adch}
and Theorem 4.1.11 \cite{henle},
since a space of rational functions $f: K\to T_yN$ such that
$f|_K$ are holomorphic is infinite-dimensional, see also
Corollaries 3.4 and 3.5 in \cite{huya}).
Using charts of the atlas of $M$ we get that ${\sf O}_{\Upsilon }(M,N)$
is infinite-dimensional.
\par In view of the Stone-Weierstrass theorem \cite{fell}
and due to the theorem about a maximum of the modulus
of a holomorphic function on a simply connected domain the space
${\sf O}_{\Upsilon }(M,Y)$ is dense in $C^0(M,Y)$, hence
${\sf O}_{\Upsilon }(M,N)$ is dense in $C^0(M,N)$.
\par {\bf 2.5. Proposition.} {\it Let $A$ and $B$ be two compact complex 
manifolds with corners. 
\par $(i)$. Then $B$ is a strong 
$C^0([0,1]\times A,A)$-retract of $A$ if and only if $B$ is a
strong $C^0([0,1], {\sf O}_{\Upsilon }(A,A))$-retract of $A$.
\par $(ii)$. Two maps $f$ and $h\in {\sf O}_{\Upsilon }(A,B)$ are
$C^0([0,1]\times A,B)$-homotopic if and only if they are
$C^0([0,1],{\sf O}_{\Upsilon }(A,B))$-homotopic.}
\par {\bf Proof.} Since $C^0([0,1],{\sf O}_{\Upsilon }(A,B))
\subset C^0([0,1]\times A,B)$, then from a strong $C^0([0,1],{\sf O}_{
\Upsilon }(A,A))$-retraction (or a  $C^0([0,1],{\sf O}_{
\Upsilon }(A,B))$-homotopy) we get a strong
$C^0([0,1]\times A,A)$-retraction (or a 
$C^0([0,1]\times A,B)$-homotopy respectively). 
\par It remains to verify the opposite implications.
\par $(ii).$ Let $F\in C^0([0,1]\times A,B).$
Consider for each $x\in [0,1]$ a compact complex submanifold
$B_{\delta }$ in $A$ such that $F(x,A)\subset B_{\delta }\subset E_{\delta }$,
where $E_{\delta (x)}=\{ z\in A:$ $d(z,F(x,A))\le \delta \} $,
$d$ denotes a metric in $A$, $0<\delta (x)<1-x$, $\delta (x)<\delta (x')$
for each $x>x'$, $d(z,Y):=inf_{y\in Y}d(z,y)$ for $Y\subset A$.
In view of Lemma 2.4 $C^0([0,1],{\sf O}_{\Upsilon }(A,{\bf C^n}))$
is dense in $C^0([0,1]\times A,{\bf C^n})$.
Then for each $0<\delta $ consider a finite atlas $At_{\delta }(M)=
\{ (U_{j,\delta },\phi _{j,\delta }) : j \} $ of $M$ such that
$\sup_j diam(U_{j,\delta })<\delta $ and
$\sup_j \sup_{(y,z\in U_{j,\delta }; |x-x'|<\delta )} \| F(x,z)-
F(x',y) \| _{C^0(A,{\bf C^n})} <\delta $, where $A$ is embedded
into $\bf C^n$. We consider restrictions $F(x,*)|_{U_{j,\delta }}$
and such $At_{\delta }(M)$ generated by $Z\in \Upsilon $.
Therefore, there exist $E_{\delta (x)}$ and mappings $H(x,*)$ such that
$H(x,A)\subset E_{\delta (x)}$ and $H(x,z)\in {\sf O}_{\Upsilon }(A,
E_{\delta (x)})$ and $\| H(x,*)-F(x,*) \| _{C^0(A,{\bf C^n})}<\epsilon (x)
<\infty $ for each $x\in [0,1]$, where $\lim_{x\to 1}\epsilon (x)=0.$
Since $f, h\in {\sf O}_{\Upsilon }(A,B)$ we can choose
$H(0,z)=f(z)$ and $H(1,z)=h(z),$
consequently, $f$ and $h$ are $C^0([0,1],{\sf O}_{\Upsilon }(
A,B))$-homotopic.
\par $(i).$ The proof of $(i)$ is analogous to that of $(ii)$
using the fact that $\bigcap_{x\in [\epsilon ,1]}E_{\delta (x)}=B$
for each $0<\epsilon <1$, since $B$ is the complex manifold 
with corners and $C^0([0,1],{\sf O}_{\Upsilon }(A,A))$ is complete.
For this let us choose a Cauchy sequence $H_n(x,z)$
in $C^0([0,1],{\sf O}_{\Upsilon }(A,A))$ instead of one $H$ as in $(ii)$
such that $H_n(0,z)=z$ for each
$z\in A$ and $H_n(1,A)=E_{\delta (1),n}$ are complex submanifolds with corners
of $A$ with $\bigcap_nE_{\delta (1),n}=B$ and 
$E_{\delta (1),n+1}\subset E_{\delta (1),n}$ for each $n\in \bf N$,
$\lim_{n\to \infty }H_n=H$
is the desired mapping.
\par {\bf 2.6. Lemma.} {\it  Let $M$ be a manifold from \S 2.1.3.3
then  there exists a mapping $q: M\to M$ such that 
\par $(i)$ $q(A_2)=\{ s_1 \} $,
\par $(ii)$ $q: (A_1\setminus A_3)\to (M\setminus \{ s_1 \} )$
is an ${\sf O}_{\Upsilon }$-diffeomorphism,
\par $(iii)$ $s_1\in A_2\setminus A_3$
\par $(iv)$ $q$ is
$C^0([0,1],{\sf O}_{\Upsilon }(M,M))$-homotopic
with the identity mapping $id_M$ on $M$.}
\par {\bf Proof.} In view of Lemma 2.5 it is sufficient to construct
$C^0([0,1]\times M,M)$-homotopy. Let $M$ be supplied with an 
atlas charts of which are diffeomorphic to subregions of
quadrants (see Lemma 2.3.1).
Then the mapping $F_2(0,z)$ from Condition $2.1.3.3(iv)$ has an extension
to $id_M$ and $F_2(1,z)$ has an extension to mapping satisfying conditions
$(i,ii,iv)$ but for $s_0$ instead of $s_1$. 
Then due to the Riemann Mapping Theorem and Central Theorem
\cite{heins} and \S \S 2.2-2.4 the mapping $F_2$ defined on
$([0,1]\times A_2)\cup (\{ 0,1 \} \times M)$ has an 
$C^0([0,1],{\sf O}_{\Upsilon }(M,M))$-extension, since $M$ is connected.
Using retraction of $A_2$
onto $\{ s_0 \} $, compositions of homotopies and 
the fact that $A_1$ and $A_2$ are simply connected we get a 
$C^0([0,1],{\sf O}_{\Upsilon }(M,M))$-homotopy of $q$ with $id_M$.
\par {\bf 2.7. Theorems. (1).}  {\it The
space $(L^MN)_{\sf O}$ from \S 2.1.4 
is the complete separable Abelian topological group.}
\par {\bf (2).} {\it It is non-trivial and non-locally compact.}
\par {\bf (3).} {\it Moreover, if 
there are two distinct points $s_0$ and $s_1$
in $A_3$, then
two groups $(L^MN)_{\sf O}$ defined for $s_0$ and $s_1$
as marked points are isomorphic.}
\par {\bf (4).} {\it $(L^MN)_{\sf O}$ is the closed proper subgroup in 
$(L^M_{\bf R}N)_{\infty ,H}$.}
\par {\bf Proof.} At first it is proved that $(S^MN)_{\sf O}$ 
is an Abelian topological monoid with the cancellation property.
\par For each $f\in {\sf O}_{\Upsilon }(M,s_0;N,y_0)$ the
range $f(M)$ is compact and connected in $N$, since $M$ is compact.
In view of Lemmas 6.8 and 6.9 \cite{swit}, \S 2.1.4 and Propositions 2.2, 2.5
and Lemma 2.4
there exists a countable subfamily $\{
Z_j:  j \in {\bf N} \} $ in $\Upsilon $ such that $Z_j\subset Z_{j+1}$ for each
$j$ and $\lim_j \tilde diam Z_j=0$. 
Therefore,
for each partition $Z$ there exists
$\delta >0$ such that for each partition $Z"\in \{ Z_j:$ $j\in {\bf N} \} $
with $\sup_i \inf_j
dist(M_i,{M"}_j) <\delta $ and for 
$f\in {\sf O}(M,s_0;N,y_0;Z)$, there exists
$f_1\in {\sf O}(M,s_0;N,y_0;Z")$, such that $fR_{\sf O}f_1$, where
$dist(A,B)=\max (\sup_{a\in A}D(a,B); \sup_{b\in B}D(b,A)),$ 
$D(a,B):=\inf_{b\in
B}d(a,b)$, $A$ and $B$ are subspaces of the metric space $\bf C^n$
with the metric $d(a,b)=|a-b|_{\bf C^n}$. 
\par If $h\in cl\{ v: fR_{\sf O}v \} $, then
there exists a net $\{ h_{\beta }\in \{ v: fR_{\sf O}v \}: \beta \in 
\alpha \}$ for an ordinal $\alpha $ such that $\lim_{\beta }h_{\beta }=h$.
In view of the preceding paragraph 
we can consider sequences instead of general nets 
and extract from a double sequence $h_{\beta ,n}$ 
a convergent to $h$ subsequence, where $\lim_nh_{\beta ,n}=h_{\beta }$ 
for each $\beta $ with $h_{\beta ,n}(x)=f_{\beta ,n}(\eta _{\beta ,n}(x))$
for each $x\in M$, $n$ and $\beta $ such that $\eta _{\beta ,n}
\in Diff^{\infty }_{s_0}(M)$,
and $f_{\beta ,n}\in \{ v: fR_{\sf O}v \} $.
Therefore, each
$R_{\sf O}$-equivalence class is closed in ${\sf O}_{\Upsilon }(M,s_0;N,y_0)$,
since ${\sf O}_{\Upsilon }(M,s_0;N,y_0)$ and $H^{\infty }_p(
M,s_0;N,y_0)$ are complete, $f(M)$ is a compact subset in $N$
and due to the definition of this equivalence relation.
Then 
\par $(i)$ $str-ind \{ {\sf O}
(M,s_0;N,y_0;Z_j); h^{Z_i}_{Z_j}; {\bf N} \} /R_{\sf O}=(L^MN)_{\sf O}$ is
separable, since each space ${\sf O}(M,s_0;N,y_0;Z_j)$ is separable
where ${\sf O}(M,s_0;N,y_0;Z):= {\sf O}_{\Upsilon }(M,s_0;N,y_0)
\cap {\sf O}(M,N;Z)$.  The
continuity of the composition follows from Notes in \S 2.1.4.
The space $str-ind \{ {\sf O}(M,s_0;N,y_0;Z_j); h^{Z_i}_{Z_j}; {\bf N} \}$ is
complete due to Theorem 12.1.4 \cite{nari}, each class 
of $R_{\sf O}$-equivalent
elements is closed in it.  
Then to each Cauchy sequence in $(S^MN)_{\sf O}$ there
corresponds a Cauchy sequence in $str-ind \{ {\sf O}_{\Upsilon }
(M,s_0;N,y_0;Z_j); h^{Z_i}_{Z_j}; {\bf N} \}$ due to
\S \S 2.2-2.5, where $Z_j$ are the corresponding
partitions of $M$.  Hence $(S^MN)_{\sf O}$ is complete.
\par In view of Lemma 2.6 mappings $f$ and $\chi ^*(f\vee w_0)$ are
$R_{\sf O}$-equivalent for each $f\in
{\sf O}_{\Upsilon }(M,s_0;N,y_0))$, where $w_0(M)=\{ y_0 \} $, since
there exists a sequence $\eta _n \in Diff^{\infty }_{s_0}(M)$
such that $\lim_{n\to \infty }diam(\eta _n(A_2))=0$ and
$w_n, f_n\in H^{\infty }_p(M,s_0;N,y_0)$ with $\lim_{n\to \infty }f_n=f$,
$\lim_{n\to \infty }w_n=w_0$
and $\lim_{n\to \infty }\chi ^*(f_n\vee w_n)(\eta _n^{-1})=f$
due to $f(s_0)=s_0$ and formulas of differentiation of composite functions
(see Theorem 2.5 in \cite{ave}).
Hence $<w_0>_{\sf O}$ is the identity element, where
$<f>_{\sf O}:=\{ h\in {\sf O}_{\Upsilon }(M,s_0;N,y_0):$ $hR_{\sf O}f \} $
denotes the equivalence class.
There exists a ${\sf O}_{\Upsilon }$-diffeomorphism $\psi : (M_1\vee M_2)\to
(M_2\vee M_1)$, where $M_j=M$ for $j=1$ and $j=2$. Therefore, 
$(h\vee f)(\psi (z))=(f\vee h)(z)$ for each $z\in (M_1\vee M_2)$ and
each $f, h\in {\sf O}_{\Upsilon }(M,s_0;N,y_0)$, consequently,
$\chi ^*(f\vee h)R_{\sf O}\chi ^* (h\vee f)$, whence
$(S^MN)_{\sf O}$ is commutative.
Evidently $\chi ^*(f\vee h)=\chi ^*(f\vee q)$ is equivalent to
$h=q$, where $f, h, q \in {\sf O}_{\Upsilon }(M,s_0;N,y_0)$.
Therefore, $<f>_{\sf O}\circ <h>_{\sf O}=<f>_{\sf O}\circ <q>_{\sf O}$
is equivalent to $<h>_{\sf O}=<q>_{\sf O}$, consequently,
$(S^MN)_{\sf O}$ has the cancellation property.
\par Let $TN$ be the tangent bundle and $\pi :  TN\to N$ be the
natural projection such that $\pi (x)=z$ for each $x\in T_zN$ with $z\in N$.
The tangent bundle $T(S^MN)_{\sf O}$ is evidently isomorphic with
$({\sf O}_{\Upsilon }(M,s_0;TN,(y_0\times 0))/R_{\sf O})\times T_{y_0}N$, 
where the equivalence relation
$R_{\sf O}$ is considered in ${\sf O}_{\Upsilon }(M,s_0;TN, (y_0\times 0))$, 
$y_0\in N$. If $f(M)\ne q(M)$, then $f$ is not equivalent to $q$.
Therefore, $T_e(S^MN)_{\sf O}$ is infinite-dimensional
topological vector space, hence $(S^MN)_{\sf O}$
is not locally compact (see Theorem (5.9.5) in \cite{nari}).
\par If there are two points $s_0$ and $s_1$ as in Theorem 2.7.3, then
the spaces ${\sf O}_{\Upsilon }(M,s_0;N,y_0)$ and
${\sf O}_{\Upsilon }(M,s_1;N,y_0)$ are isomorphic. On the other hand,
$\chi ^*(f\vee w_0)$ is in the class of $R_{\sf O}$-equivalent elements
$<f>_{\sf O}$ and $A_2$ is $C^0([0,1]\times M,M)$-contractible
into $s_0$ and also into $s_1$. Therefore, the monoid $(S^MN)_{\sf O}$
does not depent on the choice of a marked point $s$ in $A_3$.
\par For a commutative 
monoid with the cancellation property
$(S^MN)_{\sf O}$ there exists a 
commutative group $(L^MN)_{\sf O}$
equal to the Grothendieck group. This group algebraically is the quotient group
$F/\sf B$, where $F$ is a free Abelian group generated by 
$(S^MN)_{\sf O}$ and $\sf B$ is a subgroup of $F$ generated by
elements $[f+g]-[f]-[g]$, $f$ and $g\in (S^MN)_{\sf O}$,
$[f]$ denotes an element of $F$ corresponding to $f$. In view of
\S 9 \cite{langal} and \cite{swan} the natural mapping 
$\gamma : (S^MN)_{\sf O}\to (L^MN)_{\sf O}$ is injective.
We supply $F$ with a topology inherited from
the Tychonoff product topology of $(S^MN)_{\sf O}^{\bf Z}$,
where each element $z$ of $F$ is $z=\sum_fn_{f,z}[f]$, 
$n_{f,z}\in \bf Z$ for each $f\in (S^MN)_{\sf O}$,
$\sum_f|n_{f,z}|<\infty $. In particular $[nf]-n[f]\in \sf B$,
hence $(L^MN)_{\sf O}$ is the complete topological group and
$\gamma $ is the topological embedding
such that $\gamma (f+g)=\gamma (f)+ \gamma (g)$
for each $f, g \in (S^MN)_{\sf O}$, $\gamma (e)=e$,
since $(z+B)\in \gamma (S^MN)_{\sf O}$, when $n_{f,z}\ge 0$
for each $f$, so in general $z=z^{+}-z^{-}$, where
$(z^{+}+B)$ and $(z^{-}+B)\in \gamma (S^MN)_{\sf O}$.
\par The manifold ${\sf O}_{\Upsilon }(M,s_0;N,y_0)$ is the closed
submanifold of $H^{\infty }_p(M,s_0;N,y_0)$ and $R_{\sf O}=
R_{\infty ,H}|_{{\sf O}_{\Upsilon }}$, hence $<f>_{\sf O}\subset
<f>_{\infty ,H}$ for each $f\in {\sf O}_{\Upsilon }(M,s_0;N,y_0)$,
where $<f>_{t,H}:=\{ v:\mbox{ }v\in H^t_p(M,s_0;N,y_0),\mbox{ }vR_{t,H}f \} $.
Therefore, there exists embedding of the manifold $(S^MN)_{\sf O}$
into $(S^M_{\bf R}N)_{\infty ,H}$ such that $(S^MN)_{\sf O}$
is closed in $(S^M_{\bf R}N)_{\infty ,H}$, since  $(S^MN)_{\sf O}$
is complete and the uniformity in it inherited from 
$(S^M_{\bf R}N)_{\infty ,H}$ concides with the initial one.
Then $<\chi ^*(f\vee h)>_{\sf O}=<f>_{\sf O}\circ <h>_{\sf O}$
and $<\chi ^*(f\vee h)>_{t,H}=<f>_{t,H}\circ <h>_{t,H}$,
consequently, $(S^MN)_{\sf O}$ is the submonoid of the monoid
$(S^M_{\bf R}N)_{\infty ,H}$ and inevitably $(L^MN)_{\sf O}$
is the closed subgroup in $(L^M_{\bf R}N)_{\infty ,H}$.
If $f(M)\ne q(M)$, then $<f>_{\infty ,H}\ne <q>_{\infty ,H}$,
where $f, q\in H^{\infty }_p(M,s_0;N,y_0)$, 
since $\inf_{\eta \in Diff^{\infty }_{s_0}
(M_{\bf R})} \| f\circ \eta -q \|_{C^0}>0$.
If $f\in {\sf O}_{\Upsilon }(M,N)$, then $f(M)$ is a complex
submanifold with corners in $N$, since $f$ is piecewise holomorphic
\cite{moret,sanger}.
On the other hand, there exists a family of 
$f\in H^{\infty }_p(M,s_0;N,y_0)$ such that $f(M)$ is not 
a complex manifold with corners. Therefore, 
$(L^M_{\bf R}N)_{\infty ,H}\setminus (L^MN)_{\sf O}$ has the cardinality
${\sf c}:=card({\bf R})$, since $card(H^{\infty }_p(M,s_0;N,y_0))=\sf c$
and $card(S^MN)_{\infty ,H}=\sf c$.
\par {\bf 2.8. Notes and Definitions. 1.} In view of \S I.5 \cite{kob}
a complex manifold $M$ considered over $\bf R$
admits a Riemann metric $g$. Due to Theorem IV.2.2
\cite{kob} there exists the Levi-Civita connection (with vanishing torsion)
of $M_{\bf R}$. Suppose $\nu $ is a measure on 
$M$ corresponding to the Riemann volume element $w$ ( $m$-form )
$\nu (dx)=w(dx)/w(M)$. 
The Riemann metric $g$ is positive definite and $w$ is non-degenerate
and non-negative, since $M$ is orientable.
\par The Christoffel
symbols $\Gamma ^k_{i,j}$ of the Levi-Civit\'a derivation
(see \S 1.8.12 \cite{kling}) are of class $C^{\infty }$ for $M$.
Then the
equivalent uniformity in $H^t(M,N)$ for $0\le t<\infty $
is given by the following base $\{
(f,g)\in (H^t(M,N))^2:  {\| (\psi _j\circ f -\psi _j\circ g)\|"}_{
H^t(M,X)}<\epsilon $, 
where $D^{\alpha }=\partial ^{|\alpha |}/\partial (x
^1)^{\alpha ^1} ...\partial (x ^{2m})^{\alpha ^{2m}}$, $\epsilon >0$, 
${{\| (\psi
_j\circ f -\psi _j\circ g)\|"}^2}_{ H^t(M,X)}:= \sum_{|\alpha |\le t}
\int_{M}|D^{\alpha } (\psi _j\circ f(x ) 
-\psi _j\circ g(x ))|^2\nu
(dx) \} $, $j\in \Lambda _N$, $X$ is the Hilbert space  over $\bf C$
either $\bf C^n$ or $l_2({\bf C})$, 
$x$ are local normal coordinates in $M_{\bf R}$
(see also \S 2.1).  Let now $(Z_j:  j\in {\bf N})$ be the sequence from
the proof of Theorem 2.7.  We consider submanifolds $M_{i,k}$ and ${M'}_{j,k}$
for each partition $Z_k$ as in \S 2.1.3.3 
(with $Z_k$ instead of $Z$), $i\in \Xi
_{Z_k}$, $j\in \Gamma _{Z_k}$, where 
$\Xi _{Z_k}$ and $\Gamma _{Z_k}$ are finite
subsets of $\bf N$. 
We supply $H^{\gamma }(M,X;Z_k)$ with the following metric 
$\rho _{k, \gamma }(y):=[\sum_{i\in
\Xi }$ $ {{\| y|_{M_{i,k}}\|"}^2}_{\gamma ,i,k}]^{1/2}$ for $y\in H^{\gamma
}(M,X;Z_k)$ and $\rho _{k, \gamma }(y)=+\infty $ in the 
contrary case, where $\Xi =\Xi
_{Z_k}$, $\infty >t\ge \gamma \in \bf N$,
$\gamma \ge m+1$, ${\|
y|_{M_{i,k}}\|"}_{\gamma ,i,k}$ is given analogously to ${\| y\|"}_{H^{\gamma
}(M,X)}$, but with $\int_{M_{i,k}}$ instead of $\int_{M}$.  
\par Let $Z^{\gamma
}(M,X)$ be the completion of $str-ind \{ H^{\gamma }(M,X;Z_j); h^{Z_i}_{Z_j};
{\bf N} \}=:Q$ relative to the following norm 
${\| y\| '}_{\gamma }:= \inf_k\rho
_{k, \gamma }(y)$, as usually
$Z^{\infty }(M,X)=\bigcap_{\gamma \in \bf N}Z^{\gamma }(M,X)$.
Let ${\bar Y}^{\infty }(M,X):= \{ f: f\in Z^{\infty }(M,X),$
${\bar \partial }f_j|_{M_{j,k}}=0$ $\mbox{for each }k \} ,$
where $f\in Z^{\infty }(M,X)$ imples $f=\sum_jf_j$ with
$f_j\in H^{\infty }(M,X;Z_j)$ for each $j\in \bf N$.
\par For a domain $W$ in $\bf C^m$, which is a complex manifold with corners,
let $Y^{\Upsilon ,a}(W,X)$ (and $Z^{\Upsilon ,a}(W,X)$)
be a subspace of those
$f\in {\bar Y}^{\infty }(W,X)$ (or $f\in Z^{\infty }(W,X)$ respectively)
for which 
$$\| f\|_{\Upsilon ,a}:=(\sum_{j=0}^{\infty }
({\| f\|^{*}}_j)^2/(j!)^a)^{1/2}<\infty ,$$ 
where ${({\| f\|^{*}}_j)^2}:=
({\| f\| '}_j)^2-
({\| f\| '}_{j-1})^2$ for $j\ge 1$ and ${{\| f\|^*}_0}={\| f\|'}_0$, 
$0<a<\infty $.
\par Using the atlases $At(M)$ and $At(N)$ as in \S 2.1 
for $M$ and $N$ of class of smoothness $Y^{\Upsilon ,b}$ with $\infty >a\ge b$
we get the uniform space $Y^{\Upsilon ,a}(M,s_0;N,y_0)$ (and also
$Z^{\Upsilon ,a}(M,s_0;N,y_0)$)
of mappings $f:  M\to N$ with $f(s_0)=y_0$ such
that $\psi _j\circ f\in Y^{\Upsilon ,a}(M,X)$ 
(or $\psi _j\circ f\in Z^{\Upsilon ,a}(M,X)$ respectively) for each $j$,
where $\sum_{p\in \Lambda _M, j\in \Lambda _N}\| f_{p,j}-(w_0)_{p,j}
\|^2_{Y^{\Upsilon ,a}(W_{p,j},X)}<\infty $ for each 
$f\in Y^{\Upsilon ,a}(M,s_0;N,y_0)$ is satisfied with
$w_0(M)=\{ y_0 \} $, since $M$ is compact.
To each equivalence class $\{ g:  gR_{\sf O}f \} 
=:<f>_{\sf O}$ there corresponds an
equivalence class $<f>_{\Upsilon ,a}:= cl(<f>_{\sf O}\cap Y^{\Upsilon ,a}
(M,s_0;N,y_0))$
(or $<f>^{\bf R}_{\Upsilon ,a}:= cl(<f>_{\infty ,H}\cap Z^{\Upsilon ,a}
(M,s_0;N,y_0))$), where the closure is taken
in $Y^{\Upsilon ,a}(M,s_0;N,y_0)$ 
(or $Z^{\Upsilon ,a}(M,s_0;N,y_0)$ respectively). 
This generates equivalence relations $R_{\Upsilon ,a}$
and $R^{\bf R}_{\Upsilon ,a}$ respectively.
We denote the quotient spaces
$Y^{\Upsilon ,a}(M,s_0;N,y_0)/R_{\Upsilon ,a}$ 
and $Z^{\Upsilon ,a}(M,s_0;N,y_0)/R^{\bf R}_{\Upsilon ,a}$ 
by $(S^MN)_{\Upsilon ,a}$ and $(S^M_{\bf R}N)_{\Upsilon ,a}$
correspondingly. Using the A. Grothendieck construction we get
the loop groups $(L^MN)_{\Upsilon ,a}$ and $(L^M_{\bf R}N)_{\Upsilon ,a}$
respectively.
\par {\bf 2.8.2.} Let $M$ be infinite-dimensional complex 
$Y^{\xi '}$-manifold
with corners modelled on $l_2({\bf C})$ such that 
\par $(i)$ there is the sequence 
of the canonically embedded
complex submanifolds $\eta _m^{m+1}:
M_m\hookrightarrow M_{m+1}$ for each $m\in \bf N$
and to $s_{0,m}$ in $M_m$ it corresponds $s_{0,m+1}=
\eta _m^{m+1}(s_{0,m})$ in $M_{m+1}$, $dim_{\bf C}M_m=n(m)$,
$0<n(m)<n(m+1)$ for each $m\in \bf N$, $\bigcup_mM_m$ is dense in $M$;
\par $(ii)$ $M$ and $At(M)$ are foliated, that is, 
\par $(\alpha )$ $u_i\circ u_j^{-1}|_{u_j(U_i\cap U_j)}\to l_2$
are of the form: $u_i\circ u_j^{-1}((z^l:$ $l\in {\bf N}))=
(\alpha _{i,j,m}(z^1,...,z^{n(m)}), \gamma _{i,j,m}(z^l:$ $l>n(m)))$
for each $m$, when $M$ is without a boundary. 
If $\partial M\ne \emptyset $ then 
\par $(\beta )$ for each boundary component $M_0$ of $M$ 
and $U_j\cap M_0\ne \emptyset $ we have $\phi _j: U_j\cap M_0\to
H_{l,Q}$, moreover, $\partial M_m\subset \partial M$ for each $m$,
where $H_{l,Q}:= \{ z\in Q_j:$ $x^{2l-1}\ge 0 \} $,
$Q_j$ is a quadrant in $l_2$ such that $Int_{l_2}H_{l,Q}\ne \emptyset $
(the interior of $H_{l,Q}$ in $l_2$), $z^l=x^{2l-1}+ix^{2l}$,
$x^j\in \bf R$, $z^l\in \bf C$ (see also \S 2.1.2.4);
\par $(iii)$ $M$ is embedded into $l_2$ as a bounded closed subset;
\par $(iv)$ each $M_m$ satisfies conditions $2.1.3.3 (i-iv)$.
\par Let $W$ be a bounded canonical closed subset
in $l_2({\bf C})$ with a continuous piecewise $C^{\infty }$-boundary and $H_m$ 
an increasing sequence of finite-dimensional subspaces
over $\bf C$, $H_m\subset H_{m+1}$ and $dim_{\bf C}H_m=n(m)$ for
each $m\in \bf N$. Then there are spaces 
$P^{\infty }_{\Upsilon ,a}(W,X):=str-ind_mY^{\Upsilon ,a}(W_m,X)$,
where $W_m=W\cap H_m$ and $X$ is a separable Hilbert space over $\bf C$.
\par Let $Y^{\xi }(W,X)$
be the completion of $P^{\infty }_{\Upsilon ,a}(W,X)$
relative to the following norm 
$$\| f\|_{\xi }:=
[\sum_{m=1}^{\infty }{\| f|_{W_m} \|^2}_{Y^{\Upsilon ,a}(W_m,X)}/
(n(m)!)^{1+c}]^{1/2},$$
where $0<c<\infty $ and $\xi =(\Upsilon ,a,c)$.
Let $M$ and $N$ be the $Y^{\Upsilon ,a',c'}$-manifolds with 
$0<a'<a$ and $0<c'<c$.
\par If $N$ is finite-dimensional complex $Y^{\Upsilon ,a'}$-manifold,
then it is also $Y^{\Upsilon ,a',c'}$-manifold.
There exists 
the strict inductive limit of loop groups
$(L^{M_m}N)_{\Upsilon ,a}=:L^m$, since
there are natural embeddings $L^m\hookrightarrow L^{m+1}$,
such that each element $f\in Y^{\Upsilon ,a}(M_m,s_{0,m};N,y_0)$ is
considered in $Y^{\Upsilon ,a}(M_{m+1},s_{0,m+1};N,y_0)$
as independent from $(z^{n(m)+1},...,z^{n(m+1)-1})$ in the local normal
coordinates $(z^1,...,z^{n(m+1)})$ of $M_{m+1}$.
We denote it $str-ind_mL^m=:(L^MN)_{\Upsilon ,a}$
and also $str-ind_mQ^m=:Q^{\infty }_{\Upsilon ,a}(N,y_0)$,
\par $str-ind_mY^{\Upsilon ,a}(M_m;N)=:Q^{\infty }_{\Upsilon ,a}(N)$,
where $Q^m:=Y^{\Upsilon ,a}(M_m,s_{0,m};N,y_0)$.
Then with the help of charts of $At(M)$ and $At(N)$ the space
$Y^{\xi }(W,X)$
induces the uniformity $\tau $ in $Q^{\infty }_{\Upsilon ,a}(N,y_0)$
and the completion of it relative to $\tau $ we denote by
$Y^{\xi }(M,s_0;N,y_0)$, where $\xi =(\Upsilon ,a,c)$ and
$\sum_{p\in \Lambda _M, j\in \Lambda _N}\| f_{p,j}-(w_0)_{p,j}
\|^2_{Y^{\xi }(W_{p,j},X)}<\infty $ for each 
$f\in Y^{\xi }(M,s_0;N,y_0)$ is supposed to be satisfied with
$w_0(M)=\{ y_0 \} $, since each $M_m$ is compact.
Therefore,
using classes of equivalent elements from $Q^{\infty }_{\Upsilon ,a}(N,y_0)$ 
and their closures in $Y^{\xi }(M,s_0;
N,y_0)$ as in \S 2.8.1 we get
the corresponding loop monoids which are denoted 
$(S^MN)_{\xi }$.  With the help of A. Grothendieck construction we get
loop groups $(L^MN)_{\xi }$. Substituting spaces $Y^{\Upsilon ,a}$
over $\bf C$ 
onto $Z^{\Upsilon ,a}$ over $\bf R$ with respective modifications we get
spaces $Z^{\Upsilon ,a,c}(M,N)$ over $\bf R$,
loop monoids $(S^M_{\bf R}N)_{\xi }$ and groups $(L^M_{\bf R}N)_{\xi }$
for the multi-index $\xi =(\Upsilon ,a,c)$.
\par Let $exp:  {\tilde T}N\to N$ be the exponential mapping, 
where ${\tilde T}N$ is a neighbourhood of $N$ in $TN$
\cite{kling}.
\par {\bf 2.9. Theorems. (1).}  {\it The space $(L^MN)_{\xi }=:G$ 
for $\xi =(\Upsilon ,a)$ or $\xi =(\Upsilon ,a,c)$ from \S 2.8
is the complete separable Abelian topological group. 
Moreover, $G$ is the dense subgroup in 
$(L^MN)_{\sf O}$ for $\xi =(\Upsilon ,a)$;
$G$ is non-discrete non-locally compact and locally connected.}
\par {\bf (2).} {\it  The space $X^{\xi }(M,N):=T_e(L^MN)_{\xi }$  
is Hilbert for each $1\le m=dim_{\bf C}M\le \infty $.}
\par {\bf (3.)} {\it Let $N$ be a complex Hilbert 
$Y^{\xi '}$-manifold with $\infty >a>a'>0$ and $\infty >c>c'>0$ 
for $\xi '=(\Upsilon ,a')$ or $\xi '=(\Upsilon ,a',c')$ respectively,
then there exists a mapping 
$\tilde E: {\tilde T}(L^MN)_{\xi }\to (L^MN)_{\xi }$ 
defined by $\tilde E_{\eta
}(v)=exp_{\eta (s)}\circ v_{\eta }$ on a neighbourhood $V_{\eta }$ 
of the zero
section in $T_{\eta }(L^MN)_{\xi }$ and it is 
a $C^{\infty }$-mapping for
$Y^{\xi '}$-manifold $N$ by $v$ onto
a neighbourhood $W_{\eta }=W_e\circ \eta $ of $\eta \in (L^MN)_{\xi }$;
$\tilde E$ is the uniform isomorphism of 
uniform spaces $V_{\eta }$ and $W_{\eta
}$, where $s\in M$, $e$ is the unit element in $G$, $v\in V_{\eta },$
$1\le m\le \infty $.}
\par {\bf (4).} {\it $(L^MN)_{\xi }$ is the 
closed proper subgroup in $(L^M_{\bf R}N)_{\xi }$.}  
\par {\bf Proof.} There are true analogs of Propositions 2.2, 2.5 and 
Lemmas 2.3.1, 2.4 for the considered here classes of smoothness $Y^{\xi }$
defined with the help of Gevrey classes
due to Theorem VI.9 \cite{chau} and \cite{henle},
using the standard procedure of an increasing sequence 
of $C^{\infty }$-domains $W_n$ in a quadrant
$Q$ with $dim_{\bf C}Q < \infty $ such that $cl(\bigcup_nW_n)=Q$.
\par For $\xi =(\Upsilon ,a)$ 
or $\xi =(\Upsilon ,a,c)$ classes $<f>_{\xi }$
are closed due to \S \S 2.3, 2.7, 2.8 and the 
corresponding analog of Lemma 2.4 for the considered class 
of smoothness, since the uniform spaces 
$Y^{\xi }(M,s_0;N,y_0)$ are complete.
The space $T_e(L^{M_m}N)_{\sf O}$ is linear, 
where $e$ is the unit
element of the groop $(L^{M_m}N)_{\sf O}$.
Then in particular $X^{\xi }(M_m,N)$ is
the Banach space with $\| f\|_{X^{\xi }(M_m,N)}=\inf_{y\in f}{\|
y\|'}_{\xi }$, where $f=<y>_{\xi }$, $y\in Y^{\xi }(M,s_0;X,0)$.  On
the other hand, $X^{\xi }(M_m,N)$ is isomorphic with the completion of
$T_e(L^{M_m}N)_{\sf O}$ by the norm
$\| f\|_{X^{\xi }(M_m,N)}$.  Then $(\rho
_{k, \gamma }(y^1+y^2))^2+(\rho _{k, \gamma }(y^1-y^2))^2= 2[(\rho _{k,
\gamma }(y^1))^2+(\rho _{k, \gamma }(y^2))^2]$ for
each $y^1, y^2\in H^{\gamma }(M,X;Z_k)$ due to the 
choices of $\nu $ and $\rho
_{k, \gamma }$. If $y\in H^{\gamma }(M,X;Z_k)$ then $\rho _{k, \gamma }
(y)=\rho _{l, \gamma }(y)$ for each
$l>k$, since $\nu ({M'}_l)=0$ and $y\in H^{\gamma }(M,X;Z_l)$.  
For each $y^1,
y^2 \in {H^{\gamma }}_p (M,s_0;X,0)$ there exists 
$Z\in \Upsilon $ such that
$y^1, y^2 \in H^{\gamma }(M,X;Z)$.  Therefore, from Equality $(i)$ 
in \S 2.7 it
follows that ${\| f_1+f_2\|^2}_{X^{\xi }(M,N)}+ 
{\| f_1-f_2\|^2}_{X^{\xi
}(M,N)}=2[{\| f_1\|^2}_{X^{\xi }(M,N)} +{\| f_2\|^2}_{X^{\xi
}(M,N)}]$ for each $f_1, f_2 \in X^{\xi }(M,N)$. 
Then ${\| f_1+f_2\|^*}_k^2+ {\| f_1-f_2\|^*}_k^2=
2[{\| f_1\|^*}_k^2+{\| f_2\|^*}_k^2]$
for each $0\le k\in \bf Z$ and each $f_1$ and $f_2\in 
Q^{\infty }_{\Upsilon ,a}(X)$ of \S 2.8,
consequently, $\| f_1+f_2\|_{X^{\xi }(M,N)}^2+
\| f_1-f_2\|_{X^{\xi }(M,N)}^2=2[
\| f_1\|_{X^{\xi }(M,N)}^2+\| f_2\|_{X^{\xi }(M,N)}^2]$. Hence the formula
$4(f_1,f_2):= { \| f_1+f_2\|^2}_{X^{\xi }(M,N)} -
{\| f_1 - f_2\|^2}_{X^{\xi }(M,N)}
- i { \| f_1 + i f_2 \|^2}_{X^{\xi }(M,N)}
+i { \| f_1 - i f_2\|^2}_{X^{\xi }(M,N)} $
gives the scalar product $(f_1,f_2)$ in 
$X^{\xi }(M,N)$ and this is the Hilbert space over $\bf C$.
\par The spaces $Y^{\xi }(M,s_0;N,y_0)$ and 
$X^{\xi }(M,N)$ are complete, consequently, $G$ is complete. 
The space $X^{\xi }(M,N)$ is separable and $\Lambda _N \subset
\bf N$, consequently, $G$ is separable.  The composition and the inversion in
$(L^{M_m}N)_{\sf O}$ induces these operations in $G$, that are continuous 
due to \S \S 2.3, 2.7, 2.8,
consequently, $G$ is the Abelian topological group.  
\par Let $\beta :  M_m\to
N$ be a $Y^{\xi }$-mapping such that $\beta (s_0)=y_0$. 
If $C_0$
is the connected component of $y_0$ in $N$ then $\beta (M_m) \subset C_0$.
On the other hand, $N$ was supposed to be connected.
In view of Theorems about extensions of functions of different classes of
smoothness \cite{seel,touger} (see also \S 2.2)
and using completions in the described above
spaces  there exists a neighbourhood $W$ of $w_0$ such that
for each $f:  M_m\times \{ 0, 1 \}\to N$ of class $Y^{\xi }$ with
$f(M_m,0)=\{ y_0 \} $ and $f(s,1)=\beta (s)$ for each $s\in M$ there exists
its $Y^{\xi }$-extension $f:  M_m\times [0,1]\to N$, where $\{ 0, 1\} :=\{
0\} \cup \{ 1\} $, $\beta \in W$,
since there exists a neighbourhood $V_0$ of $y_0$ in $N$
such that it is $C^0([0,1]\times V_0,N)$-contractible into
a point.
Hence for each class $<\beta >_{\xi }$ 
in a sufficiently small (open) neighbourhood of $e$ there exists a
continuous curve $h:  [0,1]\to G$ such that 
$h(0)=e$ and $h(1)=<\beta >_{\xi }$.
\par By Theorem 2.7 the tangent space $T_eG$ is
infinite-dimensional over $\bf C$, 
consequently, $G$ is not locally compact, where
$e$ is the unit element in $G$.  
\par Let $\nabla $ be a covariant
differentiation in $N$ corresponding to the Levi-Civit\'a 
connection in $N$ due to
Theorem 5.1 \cite{flas}.  
This is possible, since $N$ is the Hilbert manifold
and hence has the partition of unity \cite{lang}.  
Therefore, there exists the
exponential mapping $exp:  {\tilde T}N\to N$ such that for each 
$z\in N$ there are a ball
$B(T_zN,0,r):=\{ y\in T_zN:  \| y\|_{T_zN}\le r \} $ 
and a neighbourhood $S_z$
of $z$ in $N$ for which $exp_z:  B(T_zN,0,r)\to S_z$ 
is the homeomorphism, since
$\phi $ is in the class of smoothness $Y^{\xi }$ due to Theorem IV.2.5
\cite{lang}, where $exp_zw=\phi (1)$, $\phi (q)$ is a geodesic, $\phi :
[0,1] \to N$, $d\phi (q)/dq|_{q=0}=w$, $\phi (0)=z$, 
$w\in B(T_zN,0,r)$ \cite{kling}. 
\par In view of Theorems 5.1 and 5.2 \cite{eliass} a mapping
\par $(i)$ $E:  {\tilde T}Y^{\xi }(M_m,s_0;N,y_0)\to
Y^{\xi }(M_m,s_0;N,y_0)$ 
is a local isomorphism, since
$a> b$, so $\sum_{k=1}^{\infty }k^2(k!)^{b-a}
<\infty $, 
where 
\par $(ii)$ $E_g(h):=exp_{g(s)}\circ h_{g}$, $s\in M_m$, $h\in TY^{\xi
}(M_m,s_0;N,y_0)$, $h_g\in T_gY^{\xi }(M_m,s_0;N,y_0)$, $g\in Y^{\xi
}(M_m,s_0;N,y_0) .$ 
In view of \cite{ebi,eliass} the tangent bundle $TY^{\xi }
(M_m,s_0;N,y_0)$ is isomorphic with $Y^{\xi }(M_m,s_0;TN, y_0\times \{ 0\}
)\times T_{y_0}N.$ 
Then $E$ induces $\tilde E$ with the help of factorisation by 
$R_{\xi }$ and the subsequent A. Grothendieck construction.
This mapping $\tilde E$ is of class of smoothness 
$C^{\infty }$ as follows
from equations for geodesics (see \S IV.3 \cite{lang}), since $TTN$ is the
$Y^{\xi "}$-manifold with $a>a">0$ and $c>c">0$. 
Indeed, this construction at first may be applied for
$(L^{M_m}N)_{\sf O}$ and then using the completion to $(L^MN)_{\xi }$.
\par The last statement is proved analogously to  that of Theorem 2.7
using classes of smoothness $Y^{\xi }$ and $Z^{\xi }$.
\par {\bf 2.10. Notes and Definitions.}  Let $l_{2,\epsilon }$ 
be the Hilbert space of
sequences $x=(x^j:$ $x^j\in {\bf C}; j\in {\bf N})$ such that $\| x\|
_{l_{2,\epsilon }}:=\{ \sum_{j=1}^{\infty } |x^j|^2j^{2\epsilon } \}
^{1/2}<\infty $. For $\epsilon =0$ we omit it as the index.
Suppose that in the either
$Y^{\Upsilon ,b}$-Hilbert or $Y^{\Upsilon ,b, d'}$-manifold
$N$ modelled on
$l_2$ (see \S 2.1 and \S 2.8) there exists a dense 
$Y^{\Upsilon ,b'}$- or $Y^{\Upsilon ,b',d"}$-Hilbert 
submanifold $N'$
modelled on $l_{2,\epsilon }$, where 
\par $(1)$ $\infty >a>b>b'>0$ and $\infty >c>d'$ and either
\par $(2)$  $\infty >\epsilon >1$ and $d'\ge d">0$  or
\par $(3)$ $\infty >\epsilon \ge 0$ and $d'>d">0$ correspondingly.
\par If $N$ is
finite-dimensional let $N'=N$. Evidently, each $Y^{\Upsilon ,b}$-manifold
is the complex $C^{\infty }$-manifold. Certainly we suppose,
that a class of smoothness of a manifold $N'$ is not less than 
that of $N$ and classes of smoothness of $M$ and $N$ are not less 
than that of a given loop group for it
as in \S 2.8 and of $G'$ as below.
Let $G':= (L^MN')_{\xi '}$ be a dense subgroup of $G=(L^MN)_{\xi }$
with
\par $(a)$ $\xi '=(\Upsilon ,a")$ such that 
$\infty >a">b$
for $\xi ={\sf O}$ and the $Y^{\Upsilon ,b}$-manifolds $M$ and $N$
and the $Y^{\Upsilon ,b'}$-manifold $N'$;
\par $(b)$ $\xi '=(\Upsilon ,a")$ such that
$a>a">b$ for $\xi =(\Upsilon ,a)$; 
\par $(c)$ $\xi '=(\Upsilon ,a",c")$ for $\xi =(\Upsilon ,a,c)$
and $dim_{\bf C}M=\infty $
such that $b<a"<a$ and $d'<c"<c$ and either
$(2)$ $\infty >\epsilon >1$ with $0<d"\le d'$ or 
$(3)$ $\infty >\epsilon \ge 0$ with $0<d"<d'$,
where $M$ and $N$ are $Y^{\Upsilon ,b,d'}$-manifolds, 
$N'$ is the $Y^{\Upsilon ,b',d"}$-manifold,
$1\le dim_{\bf C}M=:m<\infty $ in the cases $(a-b)$.
For the corresponding pair 
$G':=(L^M_{\bf R}N')_{\xi '}$ and $G:=(L^M_{\bf R}N)_{\xi }$
let indices in $(1-3)$ and $(a-c)$ be the same with substitution
of $\xi =\sf O$ on $\xi =(\infty ,H)$.
\par {\bf 2.11. Theorem.} {\it On the group $G$ 
there exists a probability quasi-invariant measure $\mu $ 
relative to a dense subgroup 
$G'$ (see \S 2.8 and \S 2.10).
Moreover, this measure can be chosen
$\infty $-continuously differentiable relative to $G'$.}
\par {\bf Proof.} {\bf (I).} From the conditions imposed on a
manifold $M$ it follows, that there exists a partition $Z$
of it into complex submanifolds with corners. Then there exists
$Z$ such that the covering induced
by it refines the covering of $At(M)$ by charts. 
There exist two mappings $f_1 $ (and 
$f_2$, when $\partial M\ne \emptyset $)
in ${\sf O}_{\Upsilon }(M,{\bf C})$ such that
$A_3=f_1^{-1}(0)$ (and $\partial M=f_2^{-1}(0)$ respectively).
The manifold $M$ is modelled on a Hilbert space $X$ over $\bf C$,
so for $A_1$ and $A_2$ there are quadrants $Q_1$ and $Q_2$
in $X$ and surjective mappings $\kappa _j: Q_j\to A_j$ such that 
$\kappa _j \in {\sf O}_{\Upsilon }(Q_j,A_j)$, moreover,
$\kappa _j$ are homeomorphisms of $Int(Q_j)$ onto $Int(A_j)$ 
and $\kappa _j(\partial Q_j)=\partial A_j$, where $j=1$ or $j=2$.
On the other hand, $\partial A_j=A_3\cup (A_j\cap \partial M)$.
The combination of two $\kappa _j$ and $Q_j$ produces
a quadrant $Q$ which is
${\sf O}_{\Upsilon }$-diffeomorphic to 
$Q_1\cup Q_2$ in $X$ and a submanifold
$S$ in $Q$ with $codim_{\bf R}S=1$ such that 
$S\subset \partial Q_1\cup \partial Q_2$, $Q_1\cap Q_2=\partial
Q_1\cap \partial Q_2$,
$\kappa _j(\partial Q_1\cap \partial Q_2)=\partial A_j$ with the
corresponding embeddings of $Q_j$ into $X$
and there exists a surjective 
${\sf O}_{\Upsilon }$-mapping $\kappa : Q\to M$ such that 
$\kappa : Int(Q\setminus S) \to M\setminus (A_3\cup \partial M)$
is a homeomorphism, where $j=1$ or $j=2$. 
Let $v=\kappa (\zeta )$ denote points
in $M$ for each $z\in Q$. 
\par {\bf (II).} At first we consider case $(b)$
for $G=(L^MN)_{\xi }$ and $G'=(L^MN')_{\xi '}$. For the compact manifold $M$
we can take $Q$ up to ${\sf O}_{\Upsilon }$-diffeomorphism equal
to $[0,1]^{2m}$
with $0$ corresponding to a marked point $s_0$
in $M$. Then the measure $\nu $ of a subset
$V_M:=\{ v:$ $v\in M,$ $card( \kappa ^{-1}(v))>1 \} $ is equal to zero
(see \S 2.8.1). There exists a continuous mapping 
\par $(i)$ $K:  {\sf O}_{\Upsilon }(M,N)\times
{\sf O}_{\Upsilon }(M,A^kN)$ 
$\to H^{\infty }_p(M,B^kM)$ given by the
following formula:
\par $(ii)$ $K(F,w)(v):=\int_0^{\zeta  ^1} ...
\int_0^{\zeta ^{2m}}(F^*w)$ 
$:=\int_{\bar F_{\zeta }}w$, 
where $A^kN=\bigoplus _{j+l=0}^k\Lambda ^{(j,l)}N$ and
$\Lambda ^{(j,l)}N$ denotes the space of
differential $(j,l)$-forms $w$ on $N$, \\
$w=\sum_{|I|=j,|J|=l}w_{I,J}dz^I\wedge d{\bar z}^J$,
where $dz^I=dz^{i_1}\wedge ...\wedge dz^{i_n}$ for a multi-index
$I=(i_1,...,i_n)$, $n\in \bf N$, $|I|=i_1+..+i_n$,
$0\le i_j\in \bf Z$, $w_{I,J}: M\to \bf C$,
$B^kN:=\bigoplus_{j=0}^k\Lambda ^jN$. Here
manifolds $A^kN$ and $B^kN$ are considered to be
of classes of smoothness $\sf O$
and $C^{\infty }$ respectively, 
where $\Lambda ^jN$ is for the manifold $N_{\bf R}$,
that is the manifold $N$ which is considered over $\bf R$, 
$\Lambda ^jN$ is the space of differential forms $w$
on $N_{\bf R}$ such that $w_I: N\to \bf C$, $w=\sum_Iw_Idx^I$
(see also Proposition 1.6.4.2 and \S 1.6.5 \cite{henle}).
It is correct, since $M$ satisfies
conditions $2.1.3.3(i-iv)$ and in local coordinates $v\in M$
\par $(iii)$ $(F^*w)_{l^1,...,l^{n+s}}(v^1,...,v^m)=
\sum_{I,J} w_{I,J} 
(\partial z^{i^1}/\partial v^{l^1})...$
$(\partial z^{i^n}/\partial v^{l^n})
(\partial \bar z^{j^1}/\partial v^{l^{n+1}})...$
$(\partial \bar z^{j^s}/\partial v^{l^{s+n}})$, 
where $z=(z^1,z^2,...)=F(v^1,...,v^m)$, 
since points $v\in M$ are paramertized
with the help of $(\zeta ^1,...,\zeta ^{2m})=\zeta \in Q$,
$v=\kappa (\zeta )$,
$z=(z^1,z^2,...)$ are local complex coordinates in $N$; 
we also write simply
${\bar F}_{\zeta }:= \{ F\circ \kappa (\tau ):
0\le \tau ^j \le \zeta ^j,$ $j=1,...,2m \} $ is a set,
$F \in {\sf O}_{\Upsilon }(M,N)$, 
$w\in {\sf O}_{\Upsilon }(M,A^kN)$,
$F^*$ is the pull back generated by $F$
(see also \cite{chen,getzl,huya}). 
We take this mapping $K$, when $dim_{\bf C}M\le dim_{\bf C}N$.
When $dim_{\bf C}M > dim_{\bf C}N$ we take  
$w\in {\sf O}_{\Upsilon }(M,A^k(N^s))$
instead of ${\sf O}_{\Upsilon }(M,A^kN)$, where $N^s=N_1\times ... \times N_s$
with $N_j=N$ for each $j=1,...,s$ such that $s\ge dim_{\bf C}M/dim_{\bf C}N$,
$s\in \bf N$.
A mapping $F: M\to N$ generates a mapping $F^{\otimes s}
:=(F,...,F): M\to N^s$ and the pull back $(F^{\otimes s})^*$ which is also
denoted simply by $F^*$, where $F^*w$ is a piecewise $C^{\infty }$-mapping 
for the considered here classes of mappings, $(F,...,F)$ is an $s$-tuplet.
\par From the definition of $Y^{\Upsilon ,a}$
in \S 2.8 it follows that $K$ has a
continuous extension 
$K:  Y^{\Upsilon ,a}(M,N)\times Y^{\Upsilon ,a}(M,A^k(N^s)) \to
Z^{\Upsilon ,a"}(M,B^kM)$ for each $0<a\le a''<\infty $
given by formula $(ii)$, where $M$ and $B^kM$ for $Z^{\Upsilon ,a"}$ 
are considered over $\bf R$, since $\bf C$ over $\bf R$ 
is isomorphic with $\bf R^2$. This is due to the
Lebesgue Theorem about the differentiation $d/dz$ of the Lebesgue indefinite
integral $\int_0^z f(x)dx$.  Let $\tilde K$ be defined on tangent spaces to
these with the help of the composition of the 
local diffeomorphism $E$ given by Formulas $2.9(i,ii)$
and $K$ as above. 
The tangent spaces over $\bf C$ can also be considered over $\bf R$.
Then $\tilde K$ is
continuously strongly differentiable such that $(D\tilde K(F,w)).
(\eta ,\psi )=\tilde K(\eta ,w)+ \tilde K(F,\psi )+\tilde K(F,L_{\eta }w)$, 
since $T_{w_0}Y^{\xi }(M,N)\subset TY^{\xi }(M,N)=Y^{\xi }(M,TN)$
for $w_0(M)=\{ y_0 \}$, also $T_{w_0}Y^{\xi }(M,A^k(N^s))\subset 
Y^{\xi }(M,A^k(T(N^s)))$ for $w_0(M)=\{ y_0\times 0 \}$, 
such that $y_0\times 0 \in
A^k(N^s)$, where $F,
\eta \in U_N\subset T_{w_0}Y^{\xi }(M,N)$ for
$\xi =(\Upsilon ,a)$
and $w, \psi \in U_k\subset T_{w_0}Y^{\xi
}(M,A^k(N^s))$, 
$U_N$ and $U_k$ are the corresponding neighbourhoods of zero
sections, $L_{\eta }w$ is the Lie derivative of $w$ along $\eta $, 
moreover, $\tilde K(F,w)\in T_{w_0}Z^{\Upsilon ,a}(M,B^kM)$.
\par In view of the formula of integration on manifolds
\cite{abma} : $\int_Mf^*w=\int_Mh^*w$
for each differential form $w$ on $N$ and 
maps $f, h\in {\sf O}_{\Upsilon }(M,s_0;N,y_0)$,
when $f\in <h>_{\sf O}$, since a subspace
of smooth functions
from $M$ to $N$ considered over $\bf R$ is dense in 
${\sf O}_{\Upsilon }(M,N)$ and $\int_Mf^*w$ is the continuous functional
on ${\sf O}_{\Upsilon }(M,N)$ for each given $w$, where
$f^*$ is defined $\nu $-almost everywhere on $M$, $M$ and $N$
are orientable (see Introduction). 
If $F(s_0)=y_0$ it does not imply
such restriction for $K(F,w)(v)|_{v=\kappa (\zeta ),
\zeta =(1,...,1)}$.  For each
continuous $w\ne 0$ there exists $F_0\in Y^{\xi }(M,N)$ with its closed
support $supp(F_0)$ such that $s_0\notin supp(F_0)$ and $K(F_0,w)(
v)|_{v=\kappa (\zeta ), \zeta =(1,...,1)}=(s_0;1,...,1)$,
since there exists $s\in M$, $s\ne
s_0$, such that $w(s)\ne 0$.  
Thus for each $w$ there are $F_0$ and a mapping
$\bar A_w:=\bar A:  {\bf R^d}\times 
Y^{\xi }(M,N) \to Y^{\xi }(M,N)$ such that
\par $(iv)$ $K(\bar A(c,F_0),w) (v)|_{v=\kappa (\zeta ), \zeta =(1,...,1)}=c$, 
where $c:=K(F,w) (v
)|_{v=\kappa (\zeta ), \zeta =(1,...,1)}$, 
$d=\sum_{l=0}^k{{2m}\choose l}=dim_{\bf R}B^k_{s_0}M$, 
${m\choose k}=m!/(k!(m-k)!)$ are binomial coefficients
(see also \cite{chen,getzl} and Hodge Decomposition Theorems
7.5.3, 7.5.5 in \cite{abma}).  On the other hand,  
$A^k(L^MN)_{\xi }=(L^MA^kN)_{\xi
}\times A^k_{y_0}N$ with the marked point $(y_0\times 0)$ in $A^kN$, 
in particular we consider $A^k(L^MM)_{\xi }$ as
$(L^MA^kM)_{\xi }\times A^k_{s_0}M$ 
which is an infinite-dimensional complex manifold
even for $k=m=1$ due to Theorem 2.9.(2),
since $dim_{\bf C}A^kM>m$
for each $1\le k\le m$, where traditionally $TM:=\bigcup_{x\in M}T_xM$,
such that $dim_{\bf C}TM=
2\times dim_{\bf C}M$ and $TU_i=U_i\times \bf X$ for 
$U_i\subset M$ corresponding to the chart in $M$, $\bf X$ is the Banach space
on which $M$ is modelled \cite{kling}.
\par The space $Y^{\xi }(M,N)\times Y^{\xi }(M,A^k(N^s)) $ is
isomorphic with $Y^{\xi }(M,N\times A^k(N^s))$,
hence $E^{-1}\circ K\circ E=\tilde K$ is
defined on a neighbourhood of the zero section in 
$TY^{\xi }(M,N\times A^k(N^s))$ 
into a neighbourhood of the zero section in 
$TZ^{\xi }(M,B^kM)$ for $\xi =(\Upsilon ,a)$. 
 The restriction of the latter mapping $\tilde K$ on the
corresponding neighbourhood of the zero section 
in $TY^{\xi }(M,s_0;N\times
A^k(N^s),y_0\times (y_0,0))$ and then the factorization by the 
equivalence relation
$R_{\xi }$ and the usage of A. Grothendieck construction
produces the mapping $K_1$ from the corresponding
neighbourhood of the zero section in $T(L^MN)_{\xi } \times
T(L^MA^k(N^s))_{\xi }$ into a neighbourhood of the zero section in
$T(L^M_{\bf R}B^kM)_{\xi }\times \bf R^d$, where 
$s\ge dim_{\bf C}M/dim_{\bf C}N$.
\par Therefore, using $\tilde E$ we get a mapping $K_1$ from
$T_e(L^MN)_{\xi } \times T_e(L^MA^k(N^s))_{\xi }$ into 
$T_e(L^M_{\bf R}B^kM)_{\xi }\times \bf R^d$ respectively
such that it is continuously
strongly differentiable with $(DK_1(f,<w>)).(\eta ,<\psi >)= K_1(\eta ,
<w>)+K_1(f,<\psi >)+K_1(f,L_{\eta }<w>)$,
where $f$, $\eta \in V_N\subset
T_e(L^MN)_{\xi }$, and 
$<w> , <\psi > \in V_k\subset T_e(L^MA^k(N^s))_{\xi }$,
$V_N$ and $V_k$ are the corresponding neighbourhoods of zero sections.  
In view
of the existence of the mapping $\tilde E$ in \S 2.9
for ${\tilde T}(L^MN)_{\xi }$ there exists the continuous 
mapping $\bar K:  W_e \times V_e \to {V'}_0$
induced by $\tilde E$ and $K_1$, where $W_e$ is a neighbourhood of
$e$ in $(L^MN)_{\xi }$,
$V_e$ is a neighbourhood of the zero section in $T_e(L^MA^k(N^s))_{\xi }$ 
for the unit element $e$ in $(L^MA^k(N^s))_{\xi }$,
${V'}_0$ is a neighbourhood of zero in the Hilbert space
$T_e(L^M_{\bf R}B^kM)_{\xi }\times {\bf R^d}$ over $\bf R$.
On the other hand, we can use the mapping
$\chi ^*$ from \S 2.1.4. This mapping
$\chi ^*$ induces the tangent mapping
$T\chi ^*: TY^{\xi }((M)_1\vee (M)_2,s_0;N,y_0)
\to TY^{\xi }(M,s_0;N,y_0)$ such that $\chi ^*$ is in the class
of smoothness $C^{\infty }$. Therefore, there is the linear mapping
(differential ) $D\chi ^* (h): T_hY^{\xi }
((M)_1\vee (M)_2,s_0;N,y_0)\to \bf F$ for each $h\in Y^{\xi }(
(M)_1\vee (M)_2,s_0;N,y_0)$, where $\bf F$ 
is the Hilbert space such that $T_zY^{\xi }(M,s_0;N,y_0)=
\{ z\}\times \bf F$ for each $z\in Y^{\xi }(M,s_0;N,y_0)$,
in particular for $z=\chi ^*(h)$ (see \cite{kling}). 
\par Then we define by induction the following mapping
\par $(v)$ $\Psi _{l,M,N}(f,<w>):=$ $\Psi _{1,M,M}(\Psi
_{l-1,M,N}(id_M,\bar K(f,<w>))$, where $\Psi _{1,M,N}(f,<w>)
:=\bar K(f,<w>$, $\Psi _{1,M,M}$ is defined analogously to $\Psi
_{1,M,N}$, but with $M$ over $\bf R$ instead of $N$ over $\bf C$, that is 
$\Psi _{1,M,M}:
(L^M_{\bf R}M)_{\xi }\times T_e(L^M_{\bf R}B^kM)_{\xi }\times {\bf R^{(l-1)d}}
\to T_e(L^M_{\bf R}B^kM)_{\xi }\times \bf R^{ld}$
due to Formula $(iv)$, 
$id_M:  M\to M$ is the
identity mapping, $id_M(z)=z$ for each $z\in M$.  
\par For a sequence of Hilbert spaces $P_q$ all over either $\bf C$ or $\bf R$
with
$q\in J\subset \bf N$ let $l_{2,\delta }(\{ P_q:  q\in J\}) $ $:=\{ x=(x^q:
x^q\in P_q, q\in J); \| x\|_{l_{2,\delta }( \{ P_q:  q\in J\} ) }$
$:=(\sum_{q\in J}{\| q^{\delta } x^q\|^2}_{P_q})^{1/2} <\infty \} $
be a new Hilbert space, where
$\infty >\delta \ge 0$, the index $\delta $ is omitted for $\delta =0$.  If $J$
is finite then $l_{2,\delta }(\{ P_q:  q\in J\}) $ is isomorphic with
$\bigotimes_{q\in J}P_q$.  
On the other hand,
$l_2$ and $l_2(\{ {\bf C^d}_q: q\in {\bf N}\} )$ are isomorphic
with equivalent norms,
since ${\| x\|_{l_2}}^2=\sum_{j=1}^{\infty }|x^j|^2=\sum_{n=0}^{\infty }
\sum_{j=1}^d|x^{nd+j}|^2={\| x\|^2}_{l_2(\{{\bf C^d}_q: q\in {\bf N}\}) }$
for each $x=\{ x^j: x^j\in {\bf C}, j\in {\bf N} \}\in l_2$.
The Hilbert space $l_{2,\epsilon }$
is isomorphic with $l_{2, \epsilon }(\{ {\bf C^d}_q: q\in {\bf N} \})$.
Their norms are equivalent,
since $(d)^{2\epsilon }{\| x\|^2}_{l_{2,\epsilon }(\{ {\bf C^d}_q:
q\in {\bf N} \} )}\ge {\| x\|^2}_{l_{2,\epsilon }}\ge {\| x\|^2}_{
l_{2,\epsilon }(\{ {\bf C^d}_q: q\in {\bf N} \} ) }$
for each $x\in l_{2,\epsilon }$.
We choose 
$P_q=T_e(L^M_{\bf R}B^kM)_{\xi }\times {\bf R^{ld}}$
for each $q$, where either $J=\bf N$ for infinite-dimensional $N$
or $J=\{ 1 \} $ for finite-dimensional $N$.
Let $l\ge 2$ be fixed for $\Psi _{l,M,N}$.
Let also $e_q\in l_2(\{ P_q:  q\in J\} )$ for
each $q\in J$ such that $\pi _q:  P_q\hookrightarrow l_2(\{ P_q:  q\in J\} )$
are the natural embeddings with 
$e_q\in \pi _q(P_q)$ and $\| e_q\|_{ l_2(\{ P_q:
q\in J\} ) }=1$ for each $q\in J$.  Therefore, due to formulas $(ii-v)$ there
exists a family $\{ <w^{i,q}>:$ $ i=1,...,d; q\in J \} $ $\subset
T_e(L^MA^k(N')^s)_{\Upsilon , \beta "}$
for $0<b<\beta " <a"<\infty $, where $k=2m=dim_{\bf R}M$, 
$<w^{i,q}>$ are the corresponding classes of equivalent elements,
such that the mapping 
$$(vi)\mbox{ }\Psi _l(f):=\sum_{q\in J}\sum_{i=1}^d 
\Psi _{l,M,N}(f,< w^{i,q}>) e_q\in
l_2(\{ P_q:  q\in J\} )$$ 
is injective.  Due to Theorem $2.9(3)$ about
properties of $\tilde E$ and the open mapping Theorem (14.4.1) \cite{nari} the
mapping $\Psi _l$ is the diffeomorphism of a suitable 
neighbourhood $U_e$ of the
unit element $e\in (L^MN)_{\xi }$ (considered as the manifold over $\bf R$)
onto a neighbourhood $V_0$ of $0$ in the
corresponding Hilbert subspace $K_0$ in $l_2(\{ P_q:  q\in J\}) $.  
Let the image $V_0$ of $U_e$ be supplied 
with the strongest uniformity relative to which $\Psi _l$ is
uniformly continuous, that produces the Hilbert space
$K_0=\bigcup_{j\in \bf N}jV_0$.
\par This follows from the consideration of a space
$Z^{\gamma }_q(M,X)$ for $dim_{\bf C}M<\infty $ which is defined to be
the completion of a subspace $f\in Z^{\gamma
+qm}(M,X)$
for which $D^{\alpha }f(v)|_{(\mbox{ there exists }j\mbox{
with } v^j=s^j_0)}=0$ relative to the following norm 
\par $\| f\|_{Z^{\gamma }_q
(M,X)}:=(\sum_{(\alpha =(\alpha ^1, ...,\alpha ^m); 0\le \alpha ^j\le
q\mbox{ for each } j=1,...,m)}$ ${\| D^{\alpha }f\| ^2}_{Z^{\gamma
}(M,X)})^{1/2}$,
where $0\le q\in \bf Z$.  In this class of smoothness
analogously to \S 2.8 we get spaces $Z^{\gamma }_q(M,N)$, 
such that $Z^{\gamma }_q
(M,s_0;N,y_0)= \{ f\in Z^{\gamma }_q(M,N):  f(s_0)=y_0 \} $
and $<f>^q_{\gamma }:= cl(<f>^{\bf R}_{\gamma +qm})$ with the closure in 
$Z^{\gamma }_q(M,s_0;N,y_0)$, where
$<f>^{\bf R}_{\gamma }$ are classes of equvalent elements 
in $Z^{\gamma }(M,s_0;N,y_0)$. Then we consider a loop group
$(L^MN)^q_{\gamma }$ constructed from the loop monoid
$(S^MN)^q_{\gamma }:=Y^{\gamma }_q(M,s_0;N,y_0)/R^q_{\gamma }$
and a Hilbert space $X^{\gamma }_q(M,N):=T_e(L^MN)^q_{\gamma }$, where
$R^q_{\gamma }$ is the equivalence relation generated by 
classes $<f>^q_{\gamma }$.
\par Then the
mapping $\Psi _l:  U_e\to K_0$ given by Formula $(vi)$ 
is continuously strongly differentiable.  There
exist neighbourhoods $V'\ni e$ in $G'$ and $U_{\xi }\ni e$ in $G$ 
such that
$V'\circ U_{\xi }\subset U_e$.  We consider next a mapping $S_{\phi
}(v):=\Psi _l\circ L_{\phi } \circ {\Psi _l}^{-1}(v)-$ $v$ with $v\in 
V_{\xi }$, $\phi \in V'$, where $V_{\xi }= \Psi _l(U_{\xi })$ and $L_{\phi
}(f):=\phi \circ f$ denotes an operator of the left shift
in $G$.  Then either $(\alpha )$
$S_{\phi }(V_{\xi })\subset l_{2,\epsilon
}(\{ {P'}_q:  q\in J\}) $ for each $\phi \in V'$, where
${P'}_q=T_e(L^M_{\bf R}B^kM)_{\Upsilon ,a''}\times \bf R^{ld}$ 
with $1<\epsilon <\infty $ and $0<d"\le d'$ 
for each $q\in J$; or $(\beta )$ $S_{\phi }(V_{\xi })
\subset l_{2,\epsilon }( \{ {P"}_q:$ $q\in J \} )$
for each $\phi \in V'$ with $0\le \epsilon <\infty $ and $0<d"<d'$,
where ${P"}_q={P'}_q$ and $\| f\|_{{P"}_q}=\| f\|_{{P'}_q}(q!)^{d'-d"}$ 
is the relation between norms of $f$ in ${P"}_q$ and ${P'}_q$ respectively
for each $q$.
\par For an open interval $J\subset \bf R$ and two Banach 
spaces $E$, $F$ and open subset $U$, $U\subset E$, if $f: J\times U\to
F$ is continuous and $D_zf(\phi ,z)$ is continuous in $J\times U$, if also
$\alpha $ and $\beta $ are two continuously differentable mappings of
$U$ into $J$, then $g(z):=\int_{\alpha (z)}^{\beta (z)}f(\phi ,z)d\phi $ is
continuously differentiable in $U$ such that $Dg(z)$ is the linear mapping
$Dg(z)(h)=(\int_{\alpha (z)}^{\beta (z)}D_zf(\phi ,
z)d\phi )h+$ $(D\beta (z)(h))
f(\beta (z),z)-$ $(D\alpha (z)(h))f(\alpha (z),z)$ 
\cite{abma}.
From Formulas $(ii-vi)$ it follows that $S_{\phi }(v)$ is strongly continuosly
differentiable by $v\in V_0$ for each $\phi \in V'$.  Moreover, 
$\partial S_{\phi }(v)/\partial v=\hat P_2\hat P_1$, where 
$(\hat P_2)^{\eta }$ for sufficiently large $\eta \in \bf N$
and $\hat P_1$ are operators of
trace class for each $\phi \in V'$ and $v\in V_{\xi }$ such that $\hat P_1:
K_0\to X'$, $\hat P_2:  X'\to l_{2,\epsilon }(\{ {P'}_q:  q\in J\}) $, $X'$ is
the corresponding Hilbert space, $\hat P_1(K_0)=X'$.  
It follows from the fact
that $Y^{\Upsilon , a}(M,N)\supset Y^{\Upsilon , b"}(M,N)\supset 
Y^{\Upsilon ,b"}(M,N')$ for each $0<b'<b<b"<a<\infty $
such that the corresponding operators of embeddings are of trace class,
since $\sum_{j\in \bf N}(j!)^{b-a"}j^l<\infty $ and either
$\sum_{j\in \bf N}
j^{-\epsilon }<\infty $ for $\epsilon >1$
or $\sum_{j=1}^{\infty }(j!)^{b'-b}<\infty $. 
In the case of $(L^MN)_{\Upsilon ,a}$ using
cylindrical functions (see below) the desired measure can be induced
from the corresponding Hilbert subspace
$X'$ either $(\alpha )$ in $l_{2,\epsilon -\delta "}(\{ {\tilde P}_q: q\in J\} )$
such that ${\tilde P}_q=T_e(L^M_{\bf R}M)_{\Upsilon ,a"}\times \bf R^{ld}$ 
for each $q$, where
$0<d"\le d'$ and $1<\epsilon -\delta "<\epsilon $;
or $(\beta )$ in $l_{2,\epsilon }(\{ {\bar P}_q: q\in J\} )$
with $0<d"<d'$ and $0\le \epsilon <\infty $, where ${\bar P}_q={\tilde P}_q$
such that $\| f\|_{{\bar P}_q}=\| f\|_{{\tilde P}_q}(q!)^{d'-d"}$ 
for each $f\in {\bar P}_q$ and each $q\in J$.
\par  There exists a
Gaussian quasimeasure $\lambda $ on $X'$.  It induces a Gaussian probability
measure $\nu $ on $K_0$ with the help of an operator of trace class $\hat P_1:
K_0\to X'$ \cite{dalf,sko}.  This measure induces a measure $\tilde \mu $ on
$U_{\xi }$ with the help of $\Psi _l$
such that $\tilde \mu (A)=\nu (\Psi _l(A))$ for each $A\in Bf(U_{\xi })$,
since $\nu (V_0)>0$.  The groups $G$ and $G'$ 
are separable and
metrizable, hence there exist locally finite coverings
$\{ \phi _i\circ W_i:$ $  i\in {\bf N} \} $ of $G$
and $\{ \phi _i\circ {W'}_i:$ $  i\in {\bf N} \} $ of $G'$
with $\phi _i\in G'$ such that $W_i$ are 
open subsets in $U_{\xi }$, ${W'}_i$ are open subsets in $V'$,
where $\phi _1=e$ and $W_1=U_{\xi }$, $e\in {W'}_i\subset W_i$
for each $i$ \cite{eng}, 
that is $\bigcup_{i\in \bf N}\phi _i\circ W_i=G$
and $\bigcup_{i\in \bf N}\phi _i\circ {W'}_i=G'$.
Then $\tilde \mu $ can be extended onto $G$ by the following formula
$\mu (A):=(\sum_{i=1}^{\infty }\tilde \mu (({\phi _i}^{-1}\circ A)\cap W_i)
2^{-i})/(\sum_{i=1}^{\infty }\tilde \mu (W_i)2^{-i})$ 
for each $A\in Bf(G)$.  In view of Theorems 26.1 and 26.2 \cite{sko} 
this $\mu $ can be chosen quasi-invariant on $G$ relative to $G'$. 
\par Therefore, to verify differentiability of $\mu $ it is sufficient
to consider $\mu $ on $W_1$ and with $\phi \in V"$ for open $V"$
in $G'$, $e\in V"\subset V'$, such that 
$\| D_{\phi }S_{\phi }(v)(X_{\phi })\|<{\bar c} \times \| v\|_{K_0}\times
\| X_{\phi }\|_{T_{\phi }G'}$ for each $\phi \in V"$, where
$0<{\bar c}=const <1$.
From the above construction and the formula for the quasi-invariance
factor $\rho _{\mu }(\phi ,g)$ with the help of $\Psi _l$ we get that
$\mu $ is $\infty $-continuously differentiable relative to $G'$, 
since $\Psi _l$ is also
$\infty $-continuously strongly differentiable due to Theorems 2.9.
Indeed, $\mu _{\phi }(E)=\int_E\mu _{\phi }(dg)$ $=\int_E\rho _{\mu }(\phi ,
g)\mu (dg)$, hence
$D_{\phi }\mu _{\phi }(E)(X_{\phi })=\int_E\{ D_{\phi }\rho _{\mu }(\phi ,
g)(X_{\phi }) \} \mu (dg)$, since ${D^k}_{\phi }\rho _{\mu }($ $\phi ,g)(
X_{1, \phi },...,X_{k, \phi })\in L^2(G, \mu ,{\bf R})$.
The latter is due to $\Psi _l\in C^{\infty }$
and since $\nu _{\phi }(A):=
\nu (U_{\phi ^{-1}}(A))$ is $\infty $-continuously differentiable
by $\phi \in G'$, where $A\in Af(V_{\xi }, \nu )$ and $U_{\phi }(v):=
v+S_{\phi }(v)$. This follows from the facts that $det(U'_{\phi })>0$
for each $\phi \in V''$,
this determinant and $U^{-1}_{\phi }$
and $exp \{ \sum_{j=1}^{\infty }[2(x-U^{-1}_{\phi }(x), e_j)
(x, e_j)- (x- U^{-1}_{\phi }(x), e_j)^2]/\lambda _j\}$
are $\infty $-continuously differentiable by $\phi \in G'$, 
where $e_j$ and 
$\lambda _j$ are the eigenvectors and eigenvalues of the correlation 
operator (self-adjoint nondegenerate positive definite
operator of trace class on the real separable Hilbert space)
of the Gaussian measure $\nu $ with mean $0$, $x\in K_0$.
\par {\bf (III).} Let now $G=(L^M_{\bf R}N)_{\infty ,H}$ and 
$G'=(L^M_{\bf R}N')_{\Upsilon , a"}$. Let $\infty >a>a">b$ 
and the measure
(denoted by ) $\tilde \mu $ on $(L^M_{\bf R}N)_{\Upsilon , a}=:G_1$ 
be from \S 2.11.(II) quite analogously to the complex case
substituting ${\sf O}_{\Upsilon }(M,N)$ on $H^{\infty }_p(M,N)$
and ${\sf O}_{\Upsilon }(M,A^k(N^s))$
on $H^{\infty }_p(M,B^k(N^s))$. 
Since by the corresponding analog of
Theorem $2.9(1)$ $G_1$ is the dense subgroup of $G$,
the measure
$\tilde \mu $ on $G_1$ induces the measure $\mu $ on $G$ such that
$\mu (Q)=\tilde \mu (
\tilde Q)$, where $\tilde Q:=[x \in G_1:$ $(h_1(x),...,h_s(x)) \in
R]$, $Q:=[x \in G:$ $(h_1(x),...,h_s(x))
\in R]$, $R \in Bf({\bf R^s})$, $h_i \in \{ h: G
\to {\bf R},$ $h \mbox{ are continuous} \} $, the real field
$\bf R$ is considered with the
standard norm. Indeed, the minimal $\sigma $-fields over $G_1$
and $G$ generated by such $\tilde Q$ and $Q$ coincide
with the Borel $\sigma $-fields
$Bf(G_1)$ and $Bf(G)$ respectively. If $Q_1\cap
Q_2=\emptyset $ for such $Q_j$ then $\tilde Q_1\cap \tilde Q_2=
\emptyset $, hence $\mu $ is additive and has $\sigma $-additive
extension on $Bf(G)$, since the uniformity in
$G_1$ is stronger than in $G$.
\par Then the spaces $C^0_b(G, {\bf R})$ and $C^0_b(G_1, {\bf R})$ of bounded
continuous functions $h: G\to \bf R$ and $h: G_1\to \bf R$ are separable
such that $\| h\|_{C^0_b(G, {\bf R})}:=\sup_{x\in G}|h(x)|<\infty $. 
There exists a countable family $\{ h_j: j\in {\bf N} \} =:F$ which is dense
in $C^0_b(G_1, {\bf R})$ and in $C^0_b(G, {\bf R})$, since if $f\in 
C^0_b(G, {\bf R})$ then its restriction $f|_{G_1}\in C^0_b(G_1, {\bf R})$.
The groups $G$ and $G_1$ are separable, consequently,
we can take $F$ separating points in $G$ and in $G_1$, since
each $h\in C^0_b(G, {\bf R})$ is entirely defined by its values on
a countable dense subset of $G$.
We may define the
following subsets of open $W_1$ in $G_1$ and
open $W$ in $G$, such that $W\cap G_1=:W_1$, $W(k,c;f):=[g\in W:$
$\rho (k;g,f) \le c]$ and
$W_1(k,c;f):=[g\in W_1: \rho (k;g,f) \le c]$, where $\infty >c>0$, $k\in
\bf N$, $f\in W_1$, the
mappings $\rho (k,k';g,f):=\sum_{h\in F(k,k')} |h(g)-h(f)|$;
$F(k,k'):=\{ h_j\in F: j=k,...,k' \}$
for each $k'>k$; $\rho (k;g,f):=\rho (1,k;f,g)$, 
so $W(k+1,c;f)\subset W(k,c;f)$ for
each $k \in \bf N$. Therefore, $\cap \{ W(k,1/k;f): k \in {\bf N} \}$
$=\{ f\}$, whence the least
$\sigma $-field $\sf A$ generated by the following family ${\sf V}:=$
$\{ W(k,c;f): c>0,
k \in {\bf N}, f\in W, W\subset G, W\mbox{ is open} \} $ 
is such that $\sf A$ $\supset Bf(G)$.
Moreover, $\bigcap_{k=1}^{\infty }(\bigcap_{m=1}^{\infty }(\bigcup_{n>m}$
$W(k,1/k;f_n)))=\{ f \} $ for each $f \in W$ and each sequence
$\{ f_n \} \subset W$ converging to $f$.
Hence $\mu (W(k,c;f)):=$ $\tilde \mu (W_1(k,c;f))$ for
each $c>0$, $f \in W_1$, $k \in N$.
From the Definition of $\tilde \mu $ it follows
that $\mu (\bigcap_{k=1}^{\infty }[\bigcap_{m=1}^{\infty }[\bigcup_{n>m}
W(k,1/k;f_n)]])=0$, consequently, $\mu $ is countably
additive on $Bf(G)$. 
\par Then with the help of such cylindrical subsets we get
that a derivative
$D^k_{\phi }\tilde \mu (W_1(j,1/j;f_n$ $))(X_{1, \phi }, ...,
X_{k, \phi })$ for each $\phi \in G'$ and $X_{i, \phi }\in \Xi (G')$
induces a measure on the family of cylindrical subsets of
$G$ coinciding with $D^k_{\phi }\mu (W(j,1/j;f_n))
(X_{1, \phi },...,X_{k, \phi })$ and it has the extension up to
the $\sigma $-additive measure 
on $\sf A$. Therefore, $\mu $ on $G$ is quasi-invariant and
$\infty $-continuouisly differentiable relative to $G'$;
analogously for $G=(L^MN)_{\sf O}$ and $G'=(L^MN')_{\Upsilon ,a"}$
\par {\bf (IV).} Consider now the case of $dim_{\bf C}M=\infty $
for $G=(L^MN)_{\xi }$ and $G'=(L^MN')_{\xi '}$. 
If $f\in Y^{\Upsilon ,a}(M_m,s_0;X,0)$ 
then $f$ is independent from local complex coordinates
$v ^j$ for each $j>n(m)=dim_{\bf C}M_m$, consequently,
$D ^{\alpha }f=0$ if there is $\alpha ^j>0$ for $l>j>n(m)$,
where $\alpha =(\alpha ^1,..., \alpha ^l)$.
Let $\nu _m$ denotes the measure $\nu $ defined on $M_m$ in \S 2.8, 
where $\nu _m(M_m)=1$.
Hence $f\in Y^{\Upsilon ,a,c}(X,0)$.
Let $A^{\infty }N$ be a $Y^{\Upsilon ,b,d'}$-submanifold
of $\bigoplus_{j+l=0}^{\infty }\Lambda ^{(j,l)}N$ and $B^{\infty }M$
be a $Y^{\Upsilon ,b,d'}$-submanifold of $\bigoplus_{j=0}^{\infty }B^jM$
such that there are natural embeddings 
$A^kN\hookrightarrow A^{\infty }N^{\infty }$
and $B^kM\hookrightarrow B^{\infty }M$ for each $k\in \bf N$.
\par There is a
mapping $K_{\infty }: Q^{\infty }_{\Upsilon ,a}(N)\times
Q^{\infty }_{\Upsilon ,a}(A^{\infty }{\tilde N})\to \bigoplus_{m\in \bf N}
Z^{\Upsilon ,a}(M_m,B^{2n(m)}M_m)$ such that $K_{\infty }(F,w):=
\{ K_m(F|_{M_m}, \omega |_{M_m}): m\in {\bf N } \} $, $K_m$ are defined for 
each $M_m$ by formula $(ii)$ of \S 2.11.(II),
where $\tilde N=N$ for $dim_{\bf C}N=\infty $ and $\tilde N$
is a submanifold of $N^{\infty }:=\bigotimes_{j=1}^{\infty }N_j$ 
with $N_j=N$ for each $j$ modelled on 
$l_{2,d'}(\{ S_j: j\in {\bf N} \} )$,
$S_j=T_yN$ for each $j$, $0<d'<c"$ (see \S \S 2.8 and 2.10,
about compositions of functions of such classes analogous to Gevrey see
\cite{chau,lurim2}).
There are the embeddings 
$\eta _m: M_m\hookrightarrow M$.
Then the atlases of $M$ and of each $M_m$ can be chosen
consistent.
Hence there are the embeddings $\chi _m$ of
$Z^{\Upsilon ,a}(M_m,B^{2n(m)}M_m)$ into $Z^{\Upsilon ,a}(M,B^{\infty }M)$.
We can choose $\chi _m$ such that $\chi _m(Z^{\Upsilon ,a}(M_m,B^{2n(m)}M_m))
\cap \chi _l(Z^{\Upsilon ,a}(M_l,B^{2n(l)}M_l))=\{ 0\} $ for each $n\ne l$,
since $M$ and $B^{\infty }M$ are the Hilbert manifolds.
Let $t-a<(a-a")/2$ and $q-c<(c-c")/2$.
Therefore, $K_{\infty }$ generates the following continuous operator
$K_{\sum }: Y^{\xi }(M;N)\times
Y^{\xi }(M;A^{\infty }{\tilde N})\to 
Z^{\Upsilon ,t,q}(M;B^{\infty }M)$ for $a\le t<\infty $
and $c<q<\infty $ given by the following formula
$K_{\sum }(F,w):=\sum_{m=1}^{\infty }\chi _m(K_m(F|_{M_m},w|_{M_m}))$,
since for the corresponding $\tilde K_m$ to $K_m$ and ${\tilde K}_{\sum }$
to $K_{\sum }$ on the tangent spaces
$\| {\tilde K}_{\sum}\| \le
[\sum_{m=1}^{\infty }{\| \tilde K_m\|^2}(n(m)!)^{c-q}]^{1/2}<\infty $,
where $M$ and $B^{\infty }M$ for $Z^{\Upsilon ,t,q}$ are considered
over $\bf R$.
Suppose $z=(z^m\in {\bf R^{d(m)}}: m\in {\bf N})$,
where $z^m=K_m(F|_{M_m},w|_{M_m})|_{
v=\kappa (\zeta ), \zeta =(1,...,1)}$, 
$d(m)=dim_{\bf R}B^{2n(m)}_{s_{0,m}}M_m$,
$z^m=(z_m^j:$ $j=1,...,d(m))$. Then 
$$\| z\|^2_{q}:=\sum_{m=1}^{\infty }
(\| z^m\|^2_m/(n(m)!)^{1+q})<\infty ,$$ where 
$\| z\|^2_m:=\sum_{j=1}^{d(m)}|z_m^j|^2$. The Hilbert space of such 
sequences we denote by ${\bar N}_{q}$. 
\par Then we use $\tilde E$ (or $E$)
as in \S 2.11(II,III) and the natural quotient mappings 
${\bar \zeta }_{\infty ,L}: Y^{\xi }(M,s_0;N,y_0)
\to (S^MN)_{\xi }$
corresponding to the classes of equivalent elements (see \S 2.8).
Therefore, we get from $K_{\sum }$ the continuous mapping 
$\hat K: 
T_e(L^MN)_{\xi }\times T_e(L^MA^{\infty }{\tilde N})_{
\xi }\to T_e(L^M_{\bf R}B^{\infty }M)_{\Upsilon ,t,q}
\times \bar N_{q}$, where the Hilbert space 
$\bar N_{q}$ appears from taking into account the operators 
$\bar A$ for each $M_m$ from \S 2.11.(II).
Then $\hat K$ and again $\tilde E$ (or $E$) generate the mapping
$\bar K: W_e\times V_e\to {V'}_0$, where $W_e$ is a neighbourhood
of $e$ in $G$, $V_e$ is a neighbourhood of the zero section in 
$T_e(L^MA^{\infty }{\tilde N})_{\xi }$,
${V'}_0$ is a neighbourhood of zero in the Hilbert space $H_{t,q}:=
T_e(L^M_{\bf R}
B^{\infty }M)_{\Upsilon ,t,q}\times \bar N_{q}$. Let $l_{2,\delta }(\{ 
P_q: q\in {\bf N} \} )$ be the same as in \S 2.11(II) and $P_q$ be equal to 
$H_{t,q}$ for each $q\in \bf N$. From the definition of $\bar K$ 
and the consideration of $Q^{\infty }_{\Upsilon ,\beta " }(N',y_0)$ 
dense in $Y^{\Upsilon ,\beta " ,\gamma " }(M,s_0;N',y_0)$
it follows the existence
of a family $\{ <w^{i,q,m}>: i=1,...,d(m); 
m\in {\bf N}$; $q\in {\bf N} \}
\subset T_e(L^MA^{\infty }{\tilde N}')_{\Upsilon , \beta " ,\gamma " }$ 
for $b<\beta " <a''$, $d'<\gamma " <c''$
such that the mapping  
$$\Psi _1(f):=\sum_{q=1}^{\infty }
\sum_{m\in \bf N}\sum_{i=1}^{d(m)}\bar K(f,<w^{i,q,m}>)e_q\in 
l_2(\{ P_q: q\in {\bf N} \})$$
is injective, where 
$A^{\infty }{\tilde N}'$ is a $Y^{\Upsilon ,b',d"}$-submanifold of 
$\bigoplus_{j+l=0}^{\infty }\Lambda ^{(j,l)}{\tilde N}'$ such that there are
natural embeddings $A^k{\tilde N}'\hookrightarrow A^{\infty }{\tilde N}'$ 
for each $k\in \bf N$; $w^{i,q,m}$ are independent
from local coordinates $(z^{n(m)+1},z^{n(m)+2},...)$ for each $i$, $q$ and $m$. 
Let $V_0=\Psi _1(U_e)$ 
be supplied with the strongest uniformity relative to which $\Psi _1$ is 
uniformly continuous, where $U_e$ is a neighbourhood of $e$ in $G$ such that
$U_e\subset W_e$. This gives the Hilbert space 
$K_0=\bigcup_{j\in \bf N}jV_0$. 
There are a neighbourhood $V'\ni e$ in $G'$ and $U_{\xi }\ni e$ in $G$
such that $V'\circ U_{\xi }\subset U_e$. Let $S_{\phi }(v):=
\Psi _1\circ L_{\phi }\circ \Psi _1^{-1}(v)-v$ with $v\in V_{\xi }$,
$\phi \in V'$, where $V_{\xi }=\Psi _1(U_{\xi })$. Then either
$(\alpha )$ 
$S_{\phi }(V_{\xi })\subset l_{2,\epsilon }(\{ {P'}_q: q\in {\bf N} \} )$
for each $\phi \in V'$, where ${P'}_q=H_{\beta ',\gamma '}$
with $1<\epsilon <\infty $,
$a"<\beta '<a$ and $c"<\gamma '<c$;
or $(\beta )$ $S_{\phi }(V_{\xi })\subset 
l_{2,\epsilon }(\{ {P"}_q: q\in {\bf N} \} )$ for each $\phi \in V'$, where
$a"<\beta '<a-\delta ,$ $c"<\gamma '<c-\delta $, $0\le \epsilon <\infty $,
${P"}_q={P'}_q$ and $\| f\|_{{P"}_q}=\| f\|_{{P'}_q}(q!)^{\delta }$
for each $q$, $0<\delta <\min (a-a", c-c" )/2$.
Moreover, $\partial S_{\phi }(v)/\partial v=\hat P_1\hat P_2$, 
where $\hat P_j$ are operators of trace class for each $\phi \in V'$ 
and $v\in V_{\xi }$ and $j\in \{ 1,2 \}$.  
The final part of the proof is analogous to 
that of \S 2.11.(II).
The remaining cases of $G=(L^M_{\bf R}N)_{\xi }$
and $G'=(L^M_{\bf R}N')_{\xi '}$ are analogous with substitutions
of $Y^{\xi }$ on $Z^{\xi }$ and $A^kN$ on $B^kN$ for each
$1\le k\le \infty $.
\section{Unitary representations of loop groups.}
\par {\bf 3.1. Theorem.} {\it Let $\mu $ be a
quasi-invariant relative to $G'$ measure on $(G,Bf(G))$ as in Theorem 2.11.  
Assume
also that $H:=L^2(G,\mu ,{\bf C})$ is the standard Hilbert space of equivalence
classes of square-integrable (by $\mu $) functions $f:  G\to \bf C$.  
Then there
exists a strongly continuous injective homomorphism 
$T^{\mu }:  G'\to U(H)$, where
$U(H)$ is the unitary group on $H$ in a topology induced from a Banach space
${\sf L}(H)$ of continuous linear operators 
$A: H\to H$ supplied with the operator norm.}
\par {\bf Proof.}  Let $f$ and $h$ be in $H$, their scalar product is given by
$(f,h):=\int_G \bar h(g)f(g)\mu (dg)$, where $f$ and $h:G\to \bf C$, $\bar h$
denotes complex conjugated $h$.  There exists the regular representation 
$T:=T^{\mu }:
G'\to U(H)$ defined by the following formula (see \S 2.1.1):  
$$T^{\mu }(z)f(g):=[\rho _{\mu }(z,g)]^{1/2}f(z^{-1}g).$$
For each fixed $z$ the
quasi-invariance factor $\rho _{\mu }(z,g)$ 
is continuous by $g$, hence $T(z)f(g)$ is
measurable, if $f(g)$ is measurable (relative to $Af(G,\mu )$ and $Bf({\bf
C})$).  Therefore, $(T(z)f(g),T(z)h(g))=\int_G \bar h(z^{-1}g)f(z^{-1}g) 
\rho _{\mu }(z,g)\mu (dg)=(f,h)$, 
consequently, $T$ is unitary.  From $\mu _{z'z}(dg)/\mu
(dg)=$ $\rho _{\mu }(z'z,g)=$ $\rho _{\mu }
(z,(z')^{-1}$ $g)\rho _{\mu }(z',g)=$ $[\mu _{z'z}(dg)/\mu
_{z'}(dg)] [\mu _{z'}(dg)/\mu (dg)]$ it follows that $T(z')T(z)=T(z'z)$ and
$T(id)=I$, $T(z^{-1})=T^{-1}(z)$.  
\par The embedding of $T_eG'$ into $T_eG$ is
the operator of trace class.  The measure $\mu $ on $G$ is induced by the
Gaussian measure on the corresponding separable Hilbert space
$K_0$ over $\bf R$.  In view of Theorems 26.1 and 26.2 \cite{sko} for
each $\delta >0$ and $ \{ f_1,...,f_n \} \subset H$ there exists a compact
subset $B$ in $G$ such that $\sum_{i=1}^n\int_{G\setminus B}|f_i(g)|^2\mu
(dg)<\delta ^2$.  Therefore, there exists an open neighbourhood $W'$ of $e$ in
$G'$ and an open neighbourhood $S$ of $e$ in $G$
such that $\rho _{\mu }(z,g)$ 
is continuous and bounded on $W'\times W'\circ S$, where
$S\subset W'\circ S\subset G$.  In view of this, Theorems 2.9 and 2.11
and the H\"older inequality
$\lim_{j\to \infty }\sum_{i=1}^n\| (T(z^j)-I)f_i \|_H=0$ for each sequence $\{
z^j:  j\in {\bf N} \} $ converging to $e$ in $G'$, $\lim_{j\to \infty }z^j=e$,
where $I$ is the unit operator on $H$.  Indeed, for each $v>0$ and a continuous
function $f:  G\to \bf C$ with $\| f\|_H=1$ there is an open neighbourhood $V$
of $id$ in $G'$ (in the topology of $G'$), such that $|\rho (z,g)-1|<v$ for each
$z \in V$ and each $g \in F$ for some open $F$ in $G$, $id \in F$ with $\mu
^f_z(G\setminus F)<v$ for each $z \in V$, where $\mu ^f(dg):=|f(g)|\mu (dg)$.
At first this can be done analogously for the corresponding Banach space from
which $\mu $ was induced, where $f\in \{ f_1,...,f_n \} $, $n\in \bf N$.  
\par In $H$
continuous functions $f(g)$ are dense, hence $|\int_G |f(g)-f(z^{-1}g)(
\rho _{\mu }(z,g))^{1/2}|^2 $ 
$\mu (dg)|<4v$ for each finite family $\{ f_j \}$ with
$\|f\|_H=1$ and $z \in V'=V\cap V"$, where $V"$ 
is an open neighbourhood of $id$
in $G'$ such that $\|f(g)-f(z^{-1}g)\|_H<v$ for each $z \in V"$, $0<v<1$,
consequently $T$ is strongly continuous (that is, $T$ is continuous relative to
the strong topology on $U(H)$ induced from ${\sf L}(H)$, 
see its definition in
\cite{fell} ).  
\par Moreover, $T$ is injective, since for each $g\ne id$ there
is $f \in C^0(G,{\bf C})\cap H$, such that $f(id)=0$, $f(g)=1$, and $\|
f\|_H>0$, so $T(f)\ne I$.  
\par {\bf 3.2. Note.} In general $T$ 
is not continuous relative to the norm
topology on $U(H)$, since for each $z\ne id \in G'$ and each $1/2>v>0$ there is
$f \in H$ with $\|f\|_H=1$, such that $\|f-T(z)f\|_H>v$, when $supp(f)$ is
sufficiently small with $(z\circ supp(f))\cap supp(f)=\emptyset $.  
\par {\bf 3.3. Theorem.} {\it Let $G$ be a loop group 
with a real probability quasi-invariant
measure $\mu $ relative to a dense subgroup $G'$ as in Theorem 2.11.  
Then $\mu $
may be chosen such that the associated regular unitary representation 
of $G'$ is irreducible.}  
\par {\bf Proof.}  Let a measure $\nu $ on the Hilbert space
$K_0$ be of the same type as in the proof of Theorem 2.11.
Let a $\nu $-measurable function $f:
H\to \bf C$ be such that $\nu ( \{ x\in K_0:  
f(x+y)\ne f(x) \} =0$ for each $y\in X_0$ with 
$f\in L^1(H,\nu ,{\bf C})$, where $\nu $ is quasi-invariant relative to 
shifts from a dense linear subspace $X_0$ in $K_0$.  Let
also $P_k :  l_2 \to L(k)$ be projectors 
such that $P_k(x)=x^k$ for each
$x=(\sum_{j\in {\bf N}}x^je_j )$, where $x^k:=\sum_{j=1}^k x^je_j$,
$x^k\in L(k)$, $L(k):=sp_{\bf R}(e_1,...,e_k)$, 
$sp_{\bf R}(e_j: j\in {\bf N}):=\{ y:  y\in l_2; 
y=\sum_{j=1}^nx^je_j; x^j\in {\bf R};
n\in {\bf N} \} $.  Since $K_0$
is isomorphic with $l_2$, then each finite-dimensional subspace $L(k)$
is complemented in $K_0$ \cite{nari}.
From the proof of
Proposition II.3.1 \cite{dalf} in view of the Fubini Theorem there exists a
sequence of cylindrical functions $f_k(x)=f_k(x^k)=\int_{K_0\ominus L(k)}
f(P_kx+(I-P_k)y)\nu _{I-P_k}(dy)$ which converges to $f$ in 
$L^1(K_0,\nu ,{\bf
C})$, where $\nu =\nu _{L(k)}\otimes \nu _{I-P_k}$, $\nu _{I-P_k}$ is the
measure on $K_0\ominus L(k)$.  
Each cylindrical function $f_k$ is $\nu $-almost
everywhere constant on $K_0$, 
since $L(k)\subset X_0$ for each $k\in \bf N$,
consequently, $f$ is $\nu $-almost everywhere constant on $K_0$.  
Let $A:=\Psi _l:  U_e \to V_0$ be the same as in \S 2.11.
From the construction of $G'$ and $\mu $ with the help of 
the local diffeomorphism $A$ and $\nu $ it follows
that, if a function $f\in L^1(G,\mu ,{\bf C})$ 
satisfies the following condition
$f^h(g)=f(g)$ $(mod $ $\mu )$ by $g\in G$ for each $h\in G'$, 
then $f(x)=const $
$( mod $ $\mu )$, where $f^h(g):=f(hg)$, $g\in G$.  
\par Let $f(g)=Ch_U(g)$ be
the characteristic function of a subset $U$, $U\subset G$, $U\in Af(G,\mu )$,
then $f(hg)=1 $ $\Leftrightarrow g\in h^{-1}U$.  If $f^h(g)=f(g)$ is true by
$g\in G$ $\mu $-almost everywhere, then $\mu (\{ g\in G:  f^h(g)\ne f(g) \}
)=0$, that is $\mu ( (h^{-1}U)\bigtriangleup U)=0$, consequently, 
the measure $\mu $ 
satisfies the condition $(P)$ from \S VIII.19.5 \cite{fell}, where
$A\bigtriangleup B:=(A\setminus B)\cup (B\setminus A)$ 
for each $A, B\subset G$.
For each subset $E\subset G$ the outer measure is bounded,
$\mu ^*(E)\le 1$, since $\mu
(G)=1$ and $\mu $ is non-negative \cite{boui}, consequently, 
there exists $F\in
Bf(G)$ such that $F\supset E$ and $\mu (F)=\mu ^*(E)$.  
This $F$ may be interpreted as the least upper bound in $Bf(G)$ relative to
the latter equality.
In view of the
Proposition VIII.19.5 \cite{fell} the measure $\mu $ is ergodic, that is for
each $U\in Af(G,\mu )$ and $F\in Af(G,\mu )$ with $\mu (U)
\times \mu (F)\ne 0$
there exists $h\in G'$ such that $\mu ((h\circ E)\cap F)\ne 0$.  
\par From Theorem I.1.2 \cite{dalf} 
it follows that $(G, Bf(G))$ is a Radon space,
since $G$ is separable and complete. Therefore, a class of compact subsets
approximates from below each measure $\mu ^f$, 
$\mu ^f(dg):=|f(g)|\mu (dg)$,
where $f\in L^2(G,\mu ,{\bf C})=:H$.
Due to the Egorov Theorem II.1.11 \cite{fell} for each $\epsilon >0$
and for each sequence
$f_n(g)$ converging to $f(g)$ for $\mu $-almost every $g\in G$,
when $n\to \infty $, there exists a compact subset $\sf K$
in $G$ such that $\mu (G\setminus {\sf K})<\epsilon $ and
$f_n(g)$ converges on $\sf K$ uniformly by $g\in \sf K$,
when $n\to \infty $.
In each Hilbert space $L^2({\bf R^n},\lambda ,{\bf R})$
the linear span of
functions $f(x)=exp[(b,x)-(ax,x)]$ is dense, where $b$ and $x\in 
\bf R^n$, $a$ is a real symmetric positive definite $n\times n$
matrix, $(*,*)$ is the standard scalar product in $\bf R^n$
and $\lambda $ is the Lebesgue measure on $\bf R^n$.
If a non-linear operator $U$ on $K_0$ satisfies conditions of Theorem 
26.1 \cite{sko}, then $\nu ^U(dx)/\nu (dx)=|det U'(U^{-1}(x))|
\rho _{\nu }(x-U^{-1}(x),x)$, where $\nu ^U(B):=\nu (U^{-1}B)$
for each $B\in Bf(K_0)$, $\rho _{\nu }(z,x)=exp\{ \sum_{l=1}^{\infty }
[2(z,e_l)(x,e_l)-(z,e_l)^2]/\lambda _l \}$ by Theorem 26.2 \cite{sko}, 
where $\lambda _l$ and $e_l$
are eigenvalues and eigenfunctions of the correlation operator 
$\hat P_1: K_0\to X'$
enumerated by $l\in \bf N$,
$z\in X_0$, $\rho _{\nu }(z,x):=\nu _z(dx)/\nu (dx)$, 
$\nu _z(B):=\nu (B-z)$ for each $B\in Bf(K_0)$. 
Hence in view of the Stone-Weierstrass Theorem A.8 \cite{fell}
an algebra ${\sf V}(Q)$ of finite pointwise products of 
functions from the following space 
$sp_{\bf C}\{ \psi (g):=(\rho (h,g))^{1/2}: h\in G' \}=:Q$ is dense in
$H$, since $\rho _{\mu }(e,g)=1$ for each $g\in G$
and $L_h: G\to G$ are diffeomorphisms of the manifold $G$, $L_h(g)=hg$.
\par For each $m\in \bf N$ there are $C^{\infty }$-curves
$\phi _j^b\in G'\cap W$, where $j=1,...,m$ and $b\in (-2,2):=\{ a: -2<a<2;
a\in {\bf R} \} $ is a parameter,
such that $\phi _j^b|_{b=0}=e$ and $\phi _j:=\phi _j^1$ and vectors 
$(\partial \phi _j^b/\partial b)|_{b=0}$ for $j=1,...,m$ are linearly
independent in $T_eG'$. Then the following condition
$det (\Psi (g))=0$ defines a submanifold $G_{\Psi }$
in $G$ of codimension over $\bf R$, 
\par $(i)$ $codim_{\bf R}G_{\Psi }\ge 1$, 
where $\Psi (g)$ is a matrix dependent from $g\in G$ with matrix elements
$\Psi _{l,j}(g):=D^{2l}_{\phi _j}(\rho (\phi _j,g))^{1/2}$.
If $f\in H$ is such that 
$(f(g),(\rho (\phi ,g))^{1/2})_{H}=0$
for each $\phi \in G'\cap W$, then differentials of these scalars products
by $\phi $ are zero. But ${\sf V}(Q)$ is dense in 
$H$ and in view of condition $(i)$ this means that $f=0$,
since for each $m$ there are $\phi _j\in G'\cap W$ such that
$det \Psi (g)\ne 0$ $\mu $-almost everywhere on $G$, $g\in G$.
If $\| f\|_{ H}>0$, then $\mu (supp(f))>0$,
consequently, $\mu (G'supp(f))=1$, since $G'U=G$ for each open
$U$ in $G$ and for each $\epsilon >0$ there exists an open $U$,
$U\supset supp(f)$, such that $\mu (U\setminus supp(f))<\epsilon $.
\par This means that the vector $f_0$ is cyclic, where $f_0\in H$ 
and $f_0(g)=1$ for each $g\in G$. 
From the construction of $\mu $ it follows that for each $f_{1,j}$ and 
$f_{2,j}\in H$ with $j=1,...,n$, $n\in \bf N$ and each $\epsilon >0$ there 
exists $h\in G'$ such that $|(T_hf_{1,j},f_{2,j})_H|\le 
\epsilon |(f_{1,j},f_{2,j})_H|$,
when $|(f_{1,j},f_{2,j})_H|>0$,
since $G$ is the Radon space by Theorem I.1.2 \cite{dalf}
and $G$ is not locally compact. 
This means that there is not any finite-dimensional 
$G'$-invariant subspace $H'$ in $H$ such that
$T_hH'\subset H'$ for each $h\in G'$ and $H'\ne \{ 0 \}$.
Hence if there is a $G'$-invariant closed subspace $H'\ne 0$
in $H$ it is isomorphic with the subspace
$L^2(V,\mu ,{\bf C})$, where $V\in Bf(G)$ with $\mu (V)>0$. 
\par Let ${\sf A}_G$ denotes a $*$-subalgebra of ${\sf L}(H)$
generated by the family of unitary operators 
$\{ T_h: h\in G' \} $. In view of the von Neumann
double commuter Theorem (see \S VI.24.2 \cite{fell})
${{\sf A}_G}"$ coincides with the weak and strong operator closures of
${\sf A}_G$ in ${\sf L}(H)$, where ${{\sf A}_G}'$
denotes the commuting algebra of ${\sf A}_G$ and ${{\sf A}_G}"=
({{\sf A}_G}')'$. 
\par We suppose that $\lambda $ is a probability 
Radon measure on $G'$ such that $\lambda $ has not any atoms and
$supp (\lambda )=G'$.
In view of the strong continuity of
the regular representation there exists the S. Bochner integral
$\int_GT_hf(g)\mu (dg)$ for each $f\in H$, which implies its existence 
in the weak (B. Pettis) sence. The measures $\mu $ and $\lambda $
are non-negative and bounded, hence $H\subset L^1(G,\mu ,{\bf C})$
and $L^2(G',\lambda ,{\bf C})\subset L^1(G',\lambda ,{\bf C})$
due to the Cauchy inequality. Therefore, we can apply below 
the Fubini Theorem (see \S II.16.3 \cite{fell}).
Let $f\in H$, then there exists a countable orthonormal base
$\{ f^j: j\in {\bf N} \} $ in $H\ominus {\bf C}f$. Then for each
$n\in \bf N$ the following set $B_n:=\{ q\in L^2(G',\lambda ,{\bf C} ):$
$(f^j,f)_H=\int_{G'}q(h)(f^j,T_hf_0)_H\lambda (dh)$ for $j=0,...,n \} $
is non-empty, since the vector $f_0$ is cyclic, where $f^0:=f$. 
There exists $\infty >R>\| f\|_H$ such that $B_n\cap B^R=:B^R_n$
is non-empty and weakly compact for each $n\in \bf N$, 
since $B^R$ is weakly compact, where
$B^R:=\{ q\in L^2(G',\lambda ,{\bf C} ): \| q\| \le R \} $
(see the Alaoglu-Bourbaki Theorem in \S (9.3.3) \cite{nari}).
Therefore, $B_n^R$ is a centered system of closed subsets
of $B^R$, that is, $\cap_{n=1}^mB^R_n\ne \emptyset $
for each $m\in \bf N$, hence it has a non-empty intersection, consequently,
there exists $q\in L^2(G',\lambda ,{\bf C})$ such that
$$(ii)\mbox{ }f(g)=\int_{G'}q(h)T_hf_0(g)\lambda (dh)$$ for $\mu $-a.e.
$g\in G$.
If $F\in L^{\infty }(G,\mu ,{\bf C})$, $f_1$ and $f_2\in H$,
then there exist $q_1$ and $q_2\in L^2(G',\lambda ,{\bf C})$
satisfying Equation $(ii)$. Therefore, 
$$(iii)\mbox{ }(f_1,Ff_2)_H=:c=
\int_G\int_{G'}\int_{G'}{\bar q}_1(h_1)q_2(h_2)\rho _{\mu }^{1/2}(h_1,g)
\rho _{\mu }^{1/2}(h_2,g)F(g)\lambda (dh_1)\lambda (dh_2)\mu (dg).$$
Let $\xi (h):=\int_G\int_{G'}\int_{G'}{\bar q_1}(h_1)q_2(h_2)
\rho _{\mu }^{1/2}(h_1,g) \rho _{\mu }^{1/2}(hh_2,g)
\lambda (dh_1)\lambda (dh_2) \mu (dg)$. 
Then there exists $\beta (h)\in L^2(G',\lambda ,{\bf C})$
such that 
\par $(iv)$ $\int_{G'}\beta (h)\xi (h)\lambda (dh)=c$.\\
To prove this we consider two cases. If $c=0$ it is sufficient
to take $\beta $ orthogonal to $\xi $ in $L^2(G',\lambda ,{\bf C})$. 
Each function $q\in L^2(G',\lambda ,{\bf C})$ 
can be written as $q=q^1-q^2+iq^3-iq^4$,
where $q^j(h)\ge 0$ for each $h\in G'$ and $j=1,...,4$,
hence we obtain the corresponding decomposition for $\xi $,
$\xi =\sum_{j,k}b^{j,k}\xi ^{j,k}$, where $\xi ^{j,k}$ corresponds to
$q_1^j$ and $q_2^k$, where $b^{j,k}\in \{ 1,-1,i,-i \}$. 
If $c\ne 0$ we can choose $(j_0,k_0)$ for which $\xi ^{j_0,k_0}\ne 0$
and 
\par $(v)$ $\beta $ is orthogonal to others $\xi ^{j,k}$ with 
$(j,k)\ne (j_0,k_0)$.\\ 
Otherwise, if $\xi ^{j,k}=0$ for each
$(j,k)$, then $q_l^j(h)=0$ for each $(l,j)$ and $\lambda $-a.e. $h\in G'$,
since $\xi (0)=\int_G\mu (dg)(\int_{G'}{\bar q_1}(h_1)\rho _{\mu }^{1/2}
(h_1,g)\lambda (dh_1))(\int _{G'}q_2(h_2)
\rho _{\mu }^{1/2}(h_2,g)\lambda ($ $dh_2))=0$ and this implies $c=0$, which 
is the contradiction with the assumption $c\ne 0$.
Hence there exists $\beta $ satisfying conditions $(iv, v)$.
\par Let $a(x)\in L^{\infty }(G,\mu ,{\bf C})$, $f$ and $g\in H$, 
$\beta (h)\in L^2(G',\lambda ,{\bf C})$. Since $L^2(G',\lambda ,{\bf C})$ 
is infinite-dimensional, then for each finite family of 
$a\in \{ a_1,...,a_m \} \subset L^{\infty }(G,\mu ,{\bf C})$,
$f\in \{ f_1,...,f_m \} \subset H$ there exists
$\beta (h)\in L^2(G',\lambda ,{\bf C})$, $h\in G'$, such that
$\beta $ is orthogonal to $\int_G{\bar f}_s(g)
[f_j(h^{-1}g)
(\rho _{\mu }(h,g))^{1/2}-f_j(g)]\mu (dg)$ for each $s,j=1,...,m$. Hence
each operator of multiplication on $a_j(g)$
belongs to ${{\sf A}_G}"$, since due to Formula $(iv)$
and Condition $(v)$ there exists $\beta (h)$
such that $(f_s,a_jf_l)=$ 
$\int_G\int_{G'}{\bar f}_s(g)\beta (h)(\rho _{\mu }
(h,g))^{1/2}f_l($ $h^{-1}g)
\lambda (dh) \mu (dg)$ $=\int_G\int_{G'} {\bar f}_s(g)
\beta (h)(T_hf_l(g))\lambda (dh)\mu (dg)$,
$\int_G{\bar f}_s(g)a_j(g)f_l(g)\mu (dg)$ 
$=\int_G \int_{G'}
{\bar f}_s(g)$ $\beta (h)$ $f_l(g)\lambda (dh)\mu (dg)=$
$(f_s,a_jf_l)$. 
Hence ${{\sf A}_G}"$ contains 
subalgebra of all operators of multiplication on functions from
$L^{\infty }(G,\mu ,{\bf C})$.
\par Let us remind the following. A Banach bundle $\sf B$ over 
a Hausdorff space $G'$ is a bundle $<B,\pi >$ over $G'$, together 
with operations and norms making each fiber $B_h$ ($h\in G'$)
into a Banach space such that conditions $BB(i-iv)$ are satisfied: 
$BB(i)$ $x\mapsto \| x\| $ is 
continuous on $B$ to $\bf R$; $BB(ii)$ the operation $+$ is 
continuous as a function on $\{ (x,y)\in B\times B: 
\pi (x)=\pi (y) \} $ to $B$; $BB(iii)$ for each $\lambda \in \bf C$,
the map $x\mapsto \lambda x$ 
is continuous on $B$ to $B$; $BB(iv)$ if $h\in G'$
and $\{ x^i\} $ is any net  of elements of $B$ such that $\| x^i\| 
\to 0$ and 
$\pi (x^i)\to h$ in $G'$, then $x^i\to 0_h$ in $B$, 
where $\pi : B\to G'$ is a bundle projection, 
$B_h:=\pi ^{-1}(h)$ is the fiber over $h$ (see \S II.13.4
\cite{fell}). If $G'$ is a Hausdorff topological group, then a Banach 
algebraic bundle over $G'$ is a Banach bundle ${\sf B}=<B,\pi >$ over $G'$
together with a binary operation $\bullet $ on $B$ satisfying
conditions $AB(i-v)$: 
$AB(i)$ $\pi (b\bullet c)=\pi (b)\pi (c)$ for $b$ and $c\in B$;
$AB(ii)$ for each $x$ and $y\in G'$ the product $\bullet $
is bilinear on $B_x\times B_y$ to $B_{xy}$;
$AB(iii)$ the product $\bullet $ on $B$ is associative;
$AB(iv)$ $\| b\bullet c\| \le \| b\| \times \| c\| $
($b, c\in B$); $AB(v)$ the map $\bullet $ is continuous on $B\times B$ 
to $B$ (see \S VIII.2.2 \cite{fell}). With $G'$ and a Banach algebra $\sf A$ 
the trivial Banach bundle ${\sf B}={\sf A}\times G'$ is associative, in 
particular let ${\sf A}=\bf C$ (see \S VIII.2.7 \cite{fell}).
\par The regular representation $T$ of $G'$ gives rise to a canonical regular
$H$-projection-valued measure $\bar P$:
$\bar P(W)f=Ch_Wf$, where $f\in H$, $W\in Bf(G)$, $Ch_W$ 
is the characteristic function of $W$. Therefore, $T_h\bar P(W)=\bar P
(h\circ W)T_h$ for each $h\in G'$ and $W\in Bf(G)$, since
$\rho _{\mu }(h,h^{-1}\circ g)\rho _{\mu }(h,g)=1=\rho _{\mu }(e,g)$ 
for each $(h,g)
\in G'\times G$, 
$Ch_W(h^{-1}\circ g)=Ch_{h\circ W}(g)$ and $T_h(\bar P(W)f(g))
=\rho _{\mu }(
h,g)^{1/2}\bar P(h\circ W)f(h^{-1}\circ g)$. Thus $<T,\bar P>$ is 
a system of imprimitivity for $G'$ over $G$, which is denoted 
${\sf T}^{\mu }$. This means that conditions
$SI(i-iii)$ are satisfied: $SI(i)$ $T$ is a unitary representation
of $G'$; $SI(ii)$ $\bar P$ is a regular 
$H$-projection-valued Borel measure on $G$ and 
$SI(iii)$ $T_h\bar P(W)=\bar P(h\circ W)T_h$ for all $h\in G'$ 
and $W\in Bf(G)$. 
\par For each $F\in L^{\infty }(G,\mu ,{\bf C})$ let $\bar \alpha _F$
be the operator in ${\sf L}(H)$ consisting
of multiplication by $F$: $\bar \alpha _F(f)=Ff$ for each $f\in H$. 
The map $F\to \bar \alpha _F$ is  an isometric $*$-isomorphism
of $L^{\infty }(G,\mu ,{\bf C})$ into ${\sf L}(H)$
(see \S VIII.19.2\cite{fell}). Therefore, Propositions 
VIII.19.2,5\cite{fell}
(using the approach of this particular case given above) are applicable
in our situation.
\par If $\bar p$ is a projection onto a closed ${\sf T}^{\mu }$-stable
subspace of $H$, then $\bar p$ commutes with all
$\bar P(W)$. Hence $\bar p$ commutes with multiplication by all
$F\in L^{\infty }(G,\mu ,{\bf C})$, so by \S VIII.19.2 \cite{fell}
$\bar p=\bar P(V)$, where $V\in Bf(G)$. Also $\bar p$ commutes with all
$T_h$, $h\in G'$, consequently, $(h\circ V)\setminus V$ and 
$(h^{-1}\circ V)\setminus V$ are $\mu $-null for each $h\in G'$, 
hence $\mu ((h\circ V)\bigtriangleup V)=0$ for all $h\in G'$. In view 
of ergodicity of $\mu $ and Proposition VIII.19.5 \cite{fell}
either $\mu (V)=0$ or $\mu (G\setminus V)=0$, hence
either $\bar p=0$ or $\bar p=I$, where $I$ is the unit operator.
Hence $T$ is the irreducible unitary representation.
\par {\bf 3.4. Theorem.} {\it There exists a 
bounded intertwining operator $V: L^2(G,\mu ,{\bf C})
\to L^2(G,\mu ',{\bf C})$ such that $VT^{\mu }(\psi )=T^{\mu '}
(\psi )V$ for each $\psi \in G'$ if and only if
$\mu $ and $\mu '$ are equivalent,
where $\mu $ and $\mu '$ are quasi-invaraint measures on $G$ relative to $G'$
and $T^{\mu }$ is the associated regular representation of $G'$ from Theorems
3.1 and 3.3.}
\par {\bf Proof.} If $\mu $ is equivalent with $\mu '$,
then $\mu (dg)/\mu '(dg):=\phi (g)$ is $\mu $-a.e positive,
which produces an intertwining operator $V$, which is
an isomorphism $V: L^2(G,\mu ,{\bf C})\to
L^2(G,\mu ',{\bf C})$ given by the following formula: 
$f(g)\mapsto f(g)\phi (g)$. 
\par It remains to verify the reverse implication.
In view of Theorem 3.3 
representations $T^{\mu }$ are irreducible.
It was proved in \S 3.3 that 
\par $(i)$ the weak closure of subalgebra
generated by the family $\{ T^{\mu }(h): h\in G' \} $
in the algebra of bounded linear operators ${\sf L}(H)$
contains all operators of multiplication on functions from the space
$L^{\infty }(G,\mu ,{\bf C})$, where $H:=L^2(G,\mu ,{\bf C})$.
If measures $\mu $ and $\mu '$ are singular, then 
\par $(ii)$ either 
$\sup_{(g \in G)}|\mu '(dg)/\mu (dg)|=\infty $
or $\sup_{(g \in G)}|\mu (dg)/\mu '
(dg)|=\infty $, where $\mu '(dg)/\mu (dg):=
\lim_{(\mu (B)\to 0, \infty >\mu (B)>0, g \in B)}\mu '(B)/\mu (B) \in
[0,\infty ]$, $[0,\infty ]:=([0,\infty )\cup \{ \infty \} )$,
$[0,\infty ):= \{ x: x\in {\bf R}, 0\le x \} $, $B\in Bf(G)$.
In view of the existence of the intertwining operator $V$ of $T^{\mu }$
with $T^{\mu '}$ there exists an isomorphism of Hilbert spaces
$\tau : L^2(G,\mu ,{\bf
C})\to L^2(G,\mu ',{\bf C})$, which has a continuous extension
to an isomorphism of Banach spaces
$\tau : L^{\infty }(G,\mu ,{\bf C})
\to L^{\infty }(G,\mu ',{\bf C})$ due to Condition
$(i)$. On the other hand, in view of Condition $(ii)$ there
exists a sequence $f_n\in L^2(G,\mu ,{\bf C})\cap L^{\infty }(G,\mu ,{\bf C})$
such that $0<C_1a_n\le b_n\le C_2a_n<\infty $ for each $n\in \bf N$
and there exist
$\lim_{n\to \infty }c_n<\infty $ and $\lim_{n\to \infty }d_n=\infty $,
where $C_1$ and $C_2$ are positive constants,
$a_n:=\| f_n\|_{L^2(G,\mu ,{\bf C})}$,
$b_n:=\| \tau f_n\|_{L^2(G,\mu ',{\bf C})}$,
$c_n:=\| f_n\|_{L^{\infty }(G,\mu ,{\bf C})}$,
$d_n:=\| \tau f_n\|_{L^{\infty }(G,\mu ',{\bf C})}$,
since there are sequences $\{ y_n: 0<y_n<\infty , n\in {\bf N} \} $
such that $\sum_n(y_n)^{-2}<\infty $, but 
$\sum_n(y_n)^{-1}=\infty $. This means
that singularity of $\mu $ with $\mu '$ leads to the contradiction,
consequently, $\mu $ and $\mu '$ are equivalent.
\par {\bf 3.5. Note.} It follows from \cite{dalf,sko}, that
on $K_0$ there is a family $P$ of orthogonal Gaussian measures
of cardinality $card (P)=card ({\bf R})=:\sf c$, 
which induce quasi-invariant measures on $G$
relative to $G'$ and have continuous quasi-invariance
factor on $G'\times G$. Therefore, there are
$\sf c$ non-equivalent unitary representations $T^{\mu }$ of $G'$ 
in $L^2(G,\mu ,{\bf C})$ due to Theorems 3.3 and 3.4.
\par {\bf 3.6. Theorem.} {\it On the loop group $G=(L^MN)_{\xi }$
and $G=(L^M_{\bf R}N)_{\xi }$
from \S 2.1.4 and \S 2.8 there exists a family 
of continuous characters $\{ \Xi \} $, which separate points of
$G$.}
\par {\bf Proof.} Since $N$ is either finite-dimensional or
the separable Hilbert manifold, then $N$ has a countable 
locally finite covering subordinated to the covering
of $N$ induced by the exponential mapping $exp: {\tilde T}N
\to N$ from a neighbourhood ${\tilde T}N$ of $N$ in $TN$
such that $exp_y: V_y\to W_y$ are local diffeomorphisms
of the corresponding neighbourhoods $V_y$ and $W_y$ 
of the zero section in $T_yN$ and of $y\in N$. 
Let $\lambda $ be equivalent with a Gaussian 
probability $\sigma $-additive measure either on the entire $T_yN$ or on its 
Hilbert subspace $P$. Each such $\lambda $ induces a family
of probability measures $\nu $ on $Bf(N)$ or its cylinder
subalgebra induced by the projection of $T_yN$ onto $P$, which may differ
by their supports. 
\par Let $T_yN_{\bf R}=:L$ be an infinite-dimensional separable 
Banach space over $\bf R$, so there exists a topolgical 
vector space $L^{\bf N} :=\prod_{j=1}^{\infty }L_j$, where 
$L_j=L$ for each $j\in {\bf N}$ \cite{nari}. Consider
a subspace $\Lambda ^{\infty }$ of a space of continuous 
$\infty $-multilinear functionals $w: L^{\bf N}\to \bf R$ 
such that $w(x+y)=w(x)+w(y)$, $w(\sigma x)=(-1)^{|\sigma |}w(x)$,
$w(x)=\lambda w(z)$ for each $x, y \in L^{\bf N}$, $\sigma \in S_{\infty }$
and $\lambda \in \bf R$, where $x= \{ x^j:$ $x^j\in L, j\in {\bf N} \} 
\in L^{\bf N}$, $z^j=x^j$ for each $j\ne k_0$ and $\lambda z^{k_0}=x^{k_0}$,
$S_{\infty }$ is a group of all bijections $\sigma : {\bf N}\to {\bf N}$ 
such that $card \{ j:$ $\sigma (j)\ne j \} <\aleph _0, $ 
$|\sigma |=1$ for $\sigma =\sigma _1...\sigma _n$ with odd $n\in \bf N$ and
pairwise transpositions $\sigma _l\ne I$, that is $\sigma _l(j_1)=j_2$,
$\sigma _l(j_2)=j_1$ and $\sigma _l|_{{\bf N}\setminus \{ j_1, j_2 \} }=I$
for the corresponding $j_1\ne j_2$,
$|\sigma |=2$ for even $n$ or $\sigma =I$. Then $\Lambda ^{\infty }$ 
(or $\Lambda ^j$) induces a vector bundle $\Lambda  ^{\infty }N_{\bf R}$
(or $\Lambda ^jN_{\bf R}$) on a manifold $N_{\bf R}$ of $\infty $-multilinear
skew-symmetric mappings over ${\sf F}(N_{\bf R})$ of 
$\Psi (N_{\bf R})^{\infty }$ (or $\Psi (N_{\bf R})^j$ respectively)
into ${\sf F}(N_{\bf R})$, where $\Psi (N_{\bf R})$ is a set of differentiable vector 
fields on $N_{\bf R}$ and ${\sf  F}(N_{\bf R})$ is an algebra of real-valued
$C^1$-functions on $N_{\bf R}$ (see also \cite{kling,kob}).
This $\Lambda ^{\infty }N_{\bf R}$ is the vector bundle of differential
$\infty $-forms on $N_{\bf R}$.
Then there exist a subfamily $\Lambda ^{\infty }_GN_{\bf R}$
of differential forms
$w$ on $N$ induced by the family $\{ \nu \} $. 
\par Let ${\bar B}^{\infty }N:=(\bigoplus_{0\le j \in \bf Z}
\Lambda ^jN_{\bf R})\oplus \Lambda ^{\infty }_GN_{\bf R}$ 
for $dim_{\bf R}N_{\bf R}=\infty $ and ${\bar B}^kN=\bigoplus_{j=0}^k
\Lambda ^jN_{\bf R}$ for each $k\in \bf N$.
We choose $w\in {\bar B}^kN$, where $k=\min (dim_{\bf R}N_{\bf R}, 
dim_{\bf R}M_{\bf R})$.
There exists its pull back $f^*w$ for each $f\in Y^{\xi }(M,N)$
(see \cite{kling} and \S 2.11). Let $E_j: S_j\to P$ be a family of continuous
linear operators from Banach spaces $S_j$ into
a Banach space $P$, then there exists
a continuous linear operator $E: l_{2,d'}(\{ S_j:$ $j\in {\bf N} \} )\to P$
such that $Ex=\sum_{j=1}^{\infty }E_jx^j$, where
$x=\{ x^j:$ $x^j\in S_j, j \in {\bf N } \} $ 
$\in l_{2,d'}(\{ S_j:$ $j\in {\bf N} \} )$,
$0<d'<\infty .$ If $f\in Y^{\xi }(M,N)$, then it induces the pull 
back operator $f^*w$ for each $w\in {\bar B}^k{\tilde N}$.
Therefore, there exists a pull back $f^*w$ for $\nu $ and $w$ either on 
$N^s$ or on $\tilde N$ instead of $N$
in the corresponding cases of $dim_{\bf C}M$ and $dim_{\bf C}N$
(see \S 2.11).
\par Moreover, to $f^*w$ a measure  $\mu _{w,f}$
on $M$ corresponds. 
Then $F_w(f):=\int_Mf^*w=\int_M
(f\circ \psi )^*w$ for each $f\in Y^{\xi }(M,s_0;N,y_0)$ and
$\psi \in Diff^{\infty }_{s_0}(M_{\bf R})
\cap Z^{\xi }(M_{\bf R},M_{\bf R})$ 
due to \S 26 \cite{sko} and \cite{abma}. Therefeore,
$F_{w}$ is continuous and constant on each class
$<f>_{\xi }$ due to \S 2.9 and \S 2.11. 
If $f^*w=0$ for each $w$ as above, then
$Df=0$. In view of $f(s_0)=y_0$ this implies that $f(M)=\{ y_0 \} $.
Hence for each $<f>_{\xi }\ne e$ there exists $w$ such that
$F_w(f)\ne 0$.
\par Let $\tilde \Xi : {\bf C}\to S^1$ be a continuous 
character of $\bf C$ (or $\tilde \Xi :\bf R
\to S^1$ respectively) as the additive group (see  
\cite{hew}).
Therefore, $\Xi (g):={\tilde \Xi }(F_w(f))$ 
is a continuous character on $G$
for each $g:=<f>_{\xi }\in G$. 
 
\end{document}